\documentclass{svjour3}
\usepackage[margin=1in]{geometry}

\usepackage{graphicx,subfigure}%
\usepackage{multirow,booktabs,caption}%
\usepackage{amsmath,amssymb,amsfonts}%
\usepackage{mathrsfs}%
\usepackage[title]{appendix}%
\usepackage{xcolor,hyperref}%
\usepackage{textcomp}%
\usepackage{manyfoot}%
\usepackage{booktabs}%
\usepackage{algorithm}%
\usepackage{algorithmicx}%
\usepackage{algpseudocode}%
\usepackage{listings}%

\newcommand{\diff}[0]{\textrm{d}}
\newcommand{\einf}[0]{\mathcal{E}^{\infty}}

\newcommand{\bsub}{\begin{subequations}}
\newcommand{\esub}{\end{subequations}$\!$}

\newcommand{\eps}[0]{\varepsilon}

\newcommand{\al}[1]{\textcolor{black}{#1}}
\newcommand{\hl}[1]{\textcolor{black}{#1}}

\title{The lightning method for the heat equation}

\date{\today} 
\author{
Hunter La Croix\thanks{ Dept.~of Applied and Computational Math \& Statistics, University of Notre Dame, Notre Dame, Indiana, 46656, USA; e-mail: e-mail: {\tt hlacroix@nd.edu }} \and 
Alan E. Lindsay\thanks{ Dept.~of Applied and Computational Math \& Statistics, University of Notre Dame, Notre Dame, Indiana, 46656, USA; e-mail: {\tt a.lindsay@nd.edu}}
}

\begin{document}

\maketitle



\begin{abstract}
%
%

This paper introduces a new method for solving the planar heat equation based on the lightning method. The lightning method is a recent development in the numerical solution of linear PDEs which expresses solutions using sums of polynomials and rational functions, or more generally as sums of fundamental solutions. The method is particularly well suited to handle domains with sharp corners where solution singularities are present. Boundary conditions are formed on a set of collocation points which is then solved as an overdetermined linear system. The approach of the present work is to utilize the Laplace transform to obtain a modified Helmholtz equation which is solved by an application of the lightning method. The numerical inversion of the Laplace transform is then performed by means of Talbot integration. Our validation of the method against existing results and multiple challenging test problems shows the method attains spectral accuracy with root-exponential convergence while being robust across a wide range of time intervals and adaptable to a variety of geometric scenarios.

\end{abstract}

\baselineskip=15pt




\markboth{H.~La~Croix, A.~E.~Lindsay.}{A lightning heat solver.}


\section{Introduction.}

\begin{figure}
    \centering
    \includegraphics[width=0.4\textwidth]{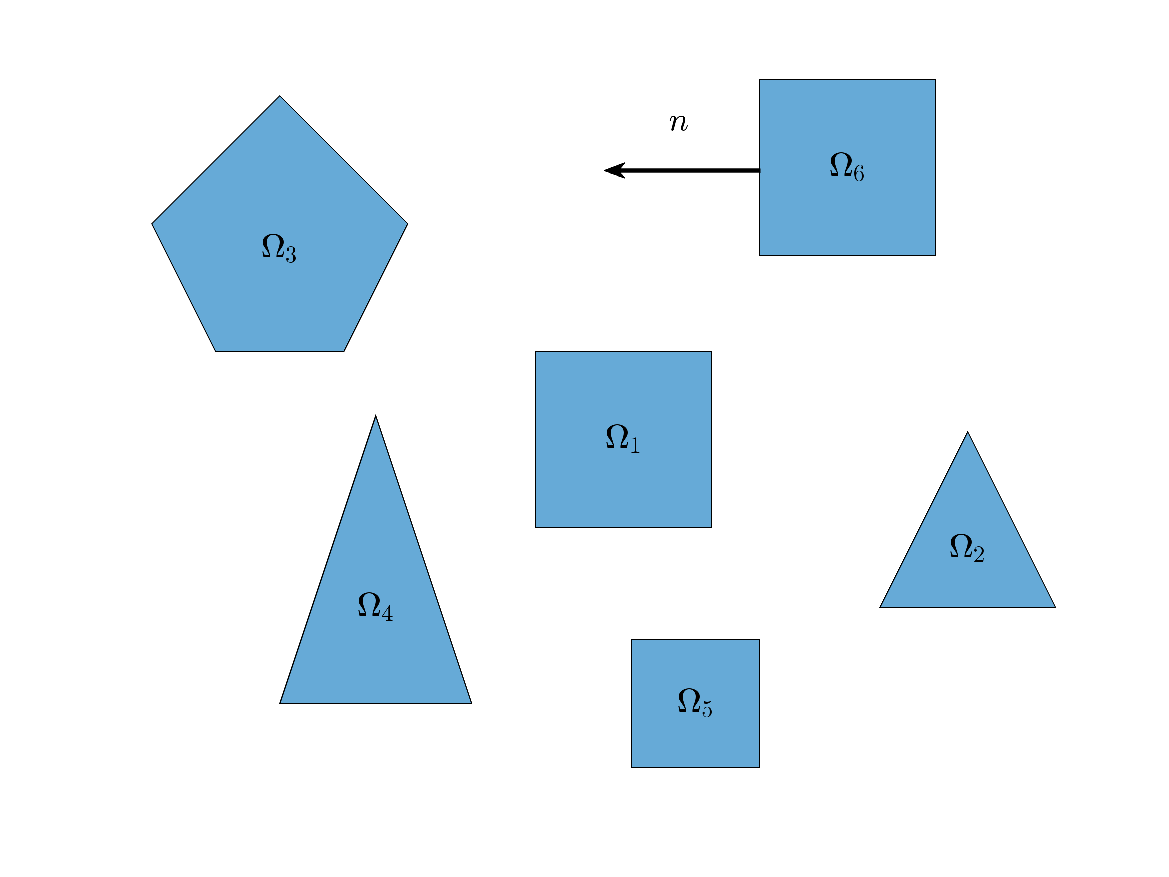}
    \caption{Schematic of the external geometry $\mathbb{R}^2\setminus\Omega$ on which the parabolic system \eqref{eq:main} is solved. The set of $N_B$ polygonal bodies is $\Omega := \cup_{k=1}^{N_B}\Omega_{k}$ and $n$ is the outward facing normal vector to $\partial\Omega$. \label{fig:Intro_Schem}}
\end{figure}

The purpose of this work is to present a new approach for the numerical solution of the planar heat equation by means of \emph{the lightning method} (LM), recently introduced by Gopal and Trefethen \cite{Gopal19A}. The specific problem of interest is the parabolic partial differential equation
\bsub\label{eq:main}
\begin{align}
    \frac{\partial u}{\partial t} &= D \Delta u, &(x,t)&\in\mathbb{R}^2\setminus\Omega \times(0,\infty);\\[4pt]
    u &= f(x),   &(x,t) & \in \partial \Omega \times(0,\infty);\\[4pt]
    u &= u_0(x), &(x,t)&\in \mathbb{R}^2 \setminus \Omega \times (t=0).
\end{align}
\esub
where $\Omega$ is a collection of $N_B$ polygonal absorbing bodies ($\Omega := \cup_{k=1}^{N_B}\Omega_{k}$). A schematic of the geometric configuration is shown in Fig.~\ref{fig:Intro_Schem}.

In this work, we show that by combining the Laplace transform with the lightning method, a rapid, robust, accurate and easy to implement numerical solution to \eqref{eq:main} can be developed \footnote{A version of our code is available at \url{https://github.com/hdlcmath/lightning-heat}}. The \al{method derived} is spectrally accurate (root-exponential convergence) and resolves corner singularities of \eqref{eq:main} inherent to polygonal geometries. On a variety of challenging test cases, we will demonstrate our method is routinely capable of attaining relative solution accuracies of order $\mathcal{O}(10^{-10})$, consistent with other implementations of the lightning method \cite{Baddoo2020,Gopal19B}.

The lightning method was originally developed in \cite{Gopal19A,Gopal19B} for the solution of Laplace and Helmholtz equations. It is in essence a complex variable technique that expresses the solution of elliptic PDEs as a rational function \al{and is hence primarily useful in two spatial dimensions}. It has since been extended to solve problems for potential flow \cite{Baddoo2020} and Stokes flow \cite{Brubeck2022}. It is able to accurately resolve solutions of linear PDEs with corners which present challenges to traditional stencil based methods due to their solution singularities.

Due to the fundamental importance of \eqref{eq:main} in many disparate applications, significant efforts have gone in to numerical techniques to approximate its solution. These including boundary integral methods \cite{cha-kre1997,fry-kro-tor2020,kro-qua2010,CHERRY2025}, Volterra integral
equations \cite{hsi-sar1993,kre1999,qiu-rie-say-zha2017} and various transform methods \cite{li-gre2007,li-gre2009}, 
Abel transforms~\cite{mel-res2014,tau2007,vee-bir2006,cos-say2004,jia-gre-wan2015}. Monte-Carlo methods are an important companion to deterministic methods that cast solutions of \eqref{eq:main} as a probability distribution that can be sampled by particle systems~\cite{MAURO2014,GROSS2022,HWANG20101089,SCHUMM2023,Lawley2024,che-lin-her-qua2022}. While slow to converge, Monte Carlo methods have no issues dealing with difficulties inherent to
stencil solution methods such as corner singularities. However, an advantage of continuous PDE methods is that the density $u(x,t)$ is obtained everywhere\al{,} which yields much more detail on the trajectories of diffusing particles. Since discrete and continuous methods for the solution of \eqref{eq:main} each have their inherent advantages and disadvantages, it is valuable to advance all feasible methods for its solution.

A particular motivation of the present study is the solution of diffusive capture problems in cellular biology \cite{Lindsay2015,Lindsay2023a,Lindsay2017Moving,che-lin-her-qua2022,LLM2020,BN2013,Lawley2024}. A specific quantity of interest is the time taken for a diffusing molecule initialized at $x_0$ to reach the target set $\partial\Omega$, the so-called \emph{first passage time} (FPT) distribution \cite{Venu2015}. This involves solving \eqref{eq:main} with  the boundary condition $f=0$ and the initial condition $u_0(x) = \delta (x-x_0)$. The main quantities of interest in this setting are
\begin{equation}\label{eqn:bdflux}
  j(t) = \int_{\partial\Omega} D \frac{\partial u}{\partial n} \, \diff s, \qquad
  c(t) = \int_{0}^{t} j(\eta) \, \diff \eta \, \al{,}
\end{equation}
 where $\partial u/\partial n := n \cdot \nabla u$ is the outward facing normal derivative. In terms of particle diffusion, this quantity is encoded by the first arrival time $\tau$ where
\bsub\label{eq:prob_dists}
\begin{equation}\label{eq:prob_dists_a}
    \tau = \inf_{0<t<\infty} \{ x(t) \in \partial \Omega \}; \qquad  \diff x = \sqrt{2D}\, dW_t, \qquad x(0) = x_0;
\end{equation}
and $W_t$ is the Wiener process. The probabilistic quantities associated \al{with} this process are then
\begin{equation}\label{eq:prob_dists_b}
    \mathbb{P}[\tau = t] = j(t), \qquad  \mathbb{P}[\tau < t] = c(t).
\end{equation}
\esub
In our validation of the numerical solution to \eqref{eq:main}, particle-based solutions  \cite{che-lin-her-qua2022} will provide a highly non-trivial but low-fidelity validation. \al{In addition, we will provide comparisons with a high order boundary integral method \cite{askham2024chunkie}.}

In the present work, we utilize the Laplace transform\al{,} which reduces \eqref{eq:main} to a modified Helmholtz equation. We then propose and implement a lighting method to solve this elliptic PDE, followed by an inversion of the Laplace transform. The inverse transform is accomplished by means of Talbot integration \cite{tre-wei-sch2006,tal1979,Weideman2015} a highly optimized and rapid method for evaluation of the Bromwich integral. This Laplace transform approach has been successfully demonstrated on parabolic problems in unbounded domains \cite{LTS2016,che-lin-her-qua2022,lindsay2024,CHERRY2025} with notable advantages being its ability to bypass the time-step limitations of traditional integrators and achieve accuracy close to machine precision.


The key enabler in this approach is a rapid and accurate method for solving the modified Helmholtz equation associated with the transform problem. In recent work \cite{CHERRY2025}, the second author and collaborators used a boundary integral approach for this purpose which expressed the solution of the modified Helmholtz equation in terms of appropriate layer potentials.  In the present work, we adapt the Lightning method for this purpose \cite{Gopal19A}. The LM expresses the solution as a series of rational functions which solve the homogeneous equation exactly. The coefficients of this series are calculated from the least squared solution of an overdetermined linear system formed at a collocation grid along $\partial \Omega$. The Lightning method has been applied to Laplace problems, whose simplicity has allowed for theoretical analysis such as proof of root-exponential convergence and optimal clustering of poles \cite{Gopal19B,Herremans23}. Previous works have demonstrated the efficacy of the LM for the Helmholtz equation \cite{Gopal19A,Ginn22}, however, its complexity when compared to the Laplace problem has hindered theoretical results. Thus, our contribution is a computational demonstration that LM solution of the Modified Helmholtz equation (and consequently the heat equation) is effective.

The outline of the paper is as follows. In Sec.~\ref{sec:num_methods} we describe our numerical approach, including implementation details of the LM for the modified Helmholtz equation and application of Talbot integration to obtain the inverse Laplace transform. In Sec.~\ref{sec:compmethods}, we will describe two alternative solutions, the kinetic Monte Carlo method and matched asymptotic expansions, \al{which} will be used to validate the results of the LM. In Sec.~\ref{sec:Results} we demonstrate our method on a variety of test problems. We demonstrate the method on domains with single or multiple absorbers, and varied boundary and initial conditions. Finally, we discuss potential extensions and avenues for future research arising from this study.

\section{The Lightning Method \& Numerical Inverse Laplace}\label{sec:num_methods}

In this section, we outline our approach to determine the numerical solution of \eqref{eq:main}. The method we describe leverages complex variables to solve linear PDEs in the plane. The solution \al{is} $u = u(z,t)$ where $z\in\mathbb{C}$, and $x= (\text{Re}(z),\text{Im}(z))$ is the independent spatial variable. The Laplace transform for $s\in\mathbb{C}$ is defined as the integral
\begin{equation}\label{eq:LaplaceTransform}
    \hat{u}(z;s) = \mathcal{L}(u) = \int_{t=0}^{\infty} u(z,t)e^{-st} \diff t.
\end{equation}
Applying \eqref{eq:LaplaceTransform} to equation \eqref{eq:main} results in the elliptic modified Helmholtz problem
\bsub\label{eq:Helmholtz}
\begin{align}
\label{eq:Helmholtz_a}    D \Delta \hat{u} - s \hat{u} &= -u_0(z), &z&\in\mathbb{C}\setminus\Omega;\\[4pt]
\label{eq:Helmholtz_b} \hat{u} &= f(z)/s,   &z& \in \partial \Omega.
\end{align}
\esub
\al{The boundary condition \eqref{eq:Helmholtz_b} arises from evaluating \eqref{eq:LaplaceTransform} at $z\in\partial\Omega$ which implies that $\hat{u}|_{\partial\Omega} =  f(z) \int_{0}^{\infty} e^{-st} \diff t = f(z)/s.$} The next step in our approach is to seek a particular solution so that the governing equation \eqref{eq:Helmholtz_a} becomes homogeneous at the expense of \al{modifying} the boundary conditions. First, we introduce the free space Green's function $G(z,\xi;s)$ satisfying
\bsub\label{eq:Greens}
\begin{gather}
	D \Delta G - s G = \hl{-}\delta(z-\xi), \qquad z \in \mathbb{C}\setminus\{\xi\};\\[5pt]
	G(z,\xi;s) =\frac{1}{2\pi D    } K_0 (\alpha|z-\xi|), \qquad \alpha = \sqrt{\frac{s}{D}}.
\end{gather}
\esub
where $K_0(z)$ is the modified Bessel function of the second kind. In terms of the Green's function \eqref{eq:Greens}, we can decompose the solution of \eqref{eq:Helmholtz} into a particular solution $\hat{u}_p$ and homogeneous solution $\hat{u}_h$. Specifically, we have that
\begin{equation} \label{eq:particularsol}
\hat{u}(z;s) = \hat{u}_p(z;s) + \hat{u}_h(z;s), \qquad \hat{u}_p(z;s) =  \int_{\mathbb{C} \setminus \Omega} G(z,\xi;s)u_0(\xi)\diff \xi.
\end{equation}
The homogeneous equation $\hat{u}_h(z;s)$ satisfies 
\bsub\label{eq:uh}
\begin{gather}
\label{eq:uh_a} D \Delta\hat{u}_h - s \hat{u}_h = 0, \qquad z \in \mathbb{C}\setminus\Omega;\\[5pt]
\label{eq:uh_b} \hat{u}_h = \tilde{f}(z;s) , \qquad z \in \partial\Omega;
\end{gather}
where the boundary term is given by
\begin{equation}\label{eq:uh_c}
    \tilde{f}(z;s) = \frac{f(z)}{s} - \int_{\mathbb{C} \setminus \Omega}G(z,\xi;s)u_0(\xi)\diff \xi\, .
\end{equation}
\esub
In the next section, we describe a method for the solution of \eqref{eq:uh}.

\subsection{The Lightning Method.} \label{sec:LM}
To begin a description on our approach, we introduce the following family of functions
\begin{eqnarray}
	\hl{\psi_k(z,\xi;s) = \frac{1}{2\pi D}K_{|k|}(\alpha|z-\xi|) \frac{(z-\xi)^k}{|z-\xi|^k}, \qquad \alpha(s) = \sqrt{\frac{s}{D}},}
\end{eqnarray}
which are exact solutions of the homogeneous modified Helmholtz equation \eqref{eq:uh_a} on $\mathbb{C} \setminus \{\xi\}$. Our method is to seek a solution to \eqref{eq:uh} in the form
\begin{equation} \label{eq:LM}
	\hl{\hat{u}_h(z;s) = \sum_{j=1}^{N_1} \big(a_j\psi_{-1}(z,z_j;s) + b_j\psi_1(z,z_j;s) \big) +
	\sum_{k=-N_2}^{N_2} c_k\psi_k(z,z_*;s)}
\end{equation}
where $a_j,b_j \in \mathbb{C}$ for $j = 1,\ldots , N_1$ and $c_k \in \mathbb{C}$ for $k= -N_2, \ldots, N_2$ are complex constants. It is important to remark that \eqref{eq:LM} is a general solution of the equation \eqref{eq:uh} and that the constant coefficients in the expansion must be chosen to enforce the boundary conditions.

The first sum is referred to as the Newman part, and is analogous to the rational terms from the lightning Laplace approximations with poles $\{ z_j\}_{j=1}^{N_1}$. The second sum is the Runge part, which can be understood as a higher order expansion about a single point $z_*$. The poles $\{ z_j\}_{j=1}^{N_1}$ are chosen to resolve singularities such as those at the corners of our geometry. The role of the Runge part is to smooth out the solution on the rest of the boundary. 

To describe the particular choice of $N_1,N_2$ for our algorithm, we first assign the $N_B$ bodies with an overall total of $N_v$ vertices/corners. We denote \al{by $m$} the number of Newman poles per corner, implying $N_1=mN_v $, and then set the order of each Runge expansion to be $N_2=\mathcal{O}(\sqrt{m})$. Due to the complex coefficients in \eqref{eq:LM}, our series has a total of $N=2N_1+2N_2+1$ terms, and $2N$ degrees of freedom.

\subsection{Choice of Newman pole clustering.}

\begin{figure}[!h]
    \centering
    \subfigure[$10$ poles per corner.]{\includegraphics[width=0.475\textwidth]{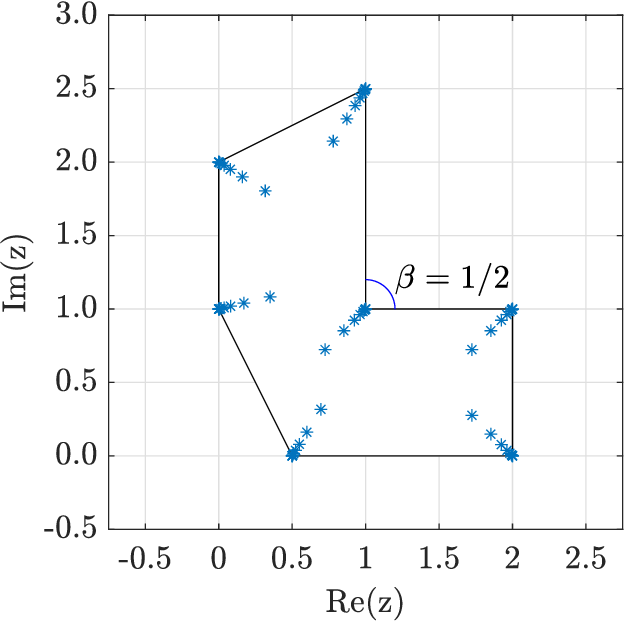}} \qquad  \subfigure[$50$ poles per corner.]{\includegraphics[width=0.475\textwidth]{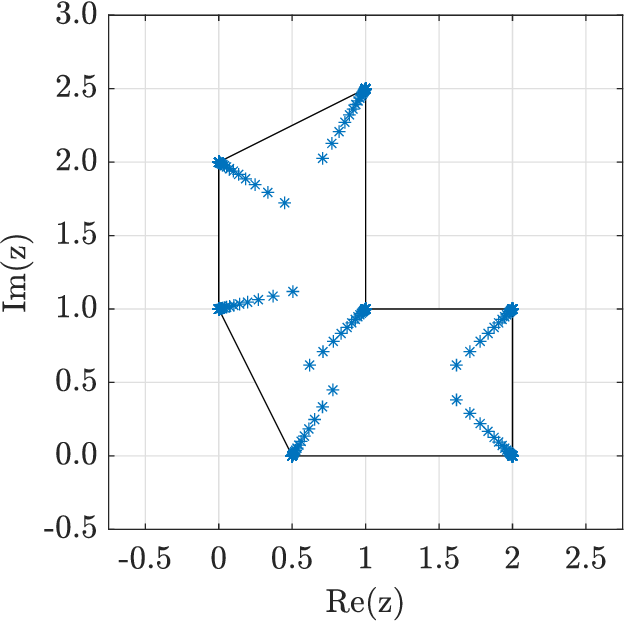}}
    \caption{Example of Newman pole clustering using the tapered distribution \eqref{eqn:clustering_dist}, with $\sigma$ in \eqref{eqn:sigma}.}\label{fig:sigma}
\end{figure}

In order to get accurate resolution of the solution and the corner singularities, one must carefully adapt the details of the clustering distributions of the Newman poles and the collocation points around them. Let us consider a corner at $z=0$, and a distribution of points that are clustering towards it from the right (positive reals). A distribution first employed by \cite{Gopal19B} was of exponential form and given by
\begin{equation}\label{eqn:clustering_dist}
    z_j = Ce^{-\sigma j/\sqrt{m}}, \qquad j=0,\ldots,m-1. 
\end{equation}
The parameter $\sigma$ controls the rate of clustering of points near the corner with empirical studies \cite{Gopal19B} suggesting the value  $\sigma =2.5$ generates a good balance between accuracy and computational effort.
A more recent analysis \hl{\cite{Trefethen2021exp}} suggested an improvement upon this, the so-called \lq\lq tapered\rq\rq\ distribution, given by
\bsub\label{eqn:tapered}
\begin{equation}\label{eqn:tapered_a}
    z_j = Ce^{-\sigma \left(\sqrt{m+1}-\sqrt{j} \right) }, \qquad j=1,\ldots,m.
\end{equation}
Moreover, for the Laplace problem, \hl{\cite{Herremans23} proposes the following geometrically-dependent choice of $\sigma$: }
\begin{equation}\label{eqn:sigma}
	\hl{\sigma = \sqrt{2(2-\beta)\beta} \pi,}
\end{equation}
\esub
where $\beta\in(0,2)$ is the exterior angle (modulo $\pi)$ of the  corner being clustered towards. \hl{This $\sigma$ is conjectured to be optimal for the Laplace problem, and we have found it to work well for the solution of the modified Helmholtz problem.}

We implemented both the clustering \al{formulas \eqref{eqn:clustering_dist} and \eqref{eqn:tapered}} for the modified Helmholtz problem \eqref{eq:uh} and found the tapered  \eqref{eqn:tapered} to perform better overall which we subsequently adopt for the reminder of this work.

For a given geometry, we choose $C=\frac{\sqrt{2}}{2} h_{\textrm{min}}$ where $h_{\textrm{min}}$ is the minimum edge length of the entire geometry $\partial\Omega$. This choice guarantees that the poles remain in the interior of the geometry. For a corner in our geometry, we find the angular bisector and move this cluster of poles to \al{lie along} it. In Fig.~\ref{fig:sigma} we show examples for $m=10$ and $m=50$ points per corner on a polygonal domain. We remark that due to floating point precision, some poles will round exactly to the vertex as $m$ increases and \hl{cause rank deficiency in our system}. To ameliorate this, we remove any poles whose distance to the vertex lies within twice machine \al{epsilon} $\varepsilon_M$. These removal events also serve as a flag that increasing $m$ further will not lead to improvements in accuracy.

\subsection{\hl{Collocation points}}

In order to form a linear system for the coefficients of the series \eqref{eq:LM}, collocation points must cluster towards corners with at least the same rate as the Newman poles. The study of \cite{Ginn22} used exponentially clustering of collocation points with a distribution $t\in[0,1]$ between the \hl{corner and midpoint} given by
\begin{equation} \label{eqn:oldcolpts}
	\hl{f(t)=t^Ae^{-\rho_c(1-t)}},
\end{equation}
with parameter values $\rho_c=4$ and varied $A$. An earlier version of the LM \cite{Gopal19A} used $A=0$ and $\rho_c=-4\sqrt{m}$. While the details of the parameters may vary, the broad idea is to create a balance of points towards the corner and midpoint of \al{the edges} such that the distribution near the corner matches that of the Newman pole distribution while not under-resolving the center.

\hl{To this end, our strategy of choosing collocation points will be using a simple exponential ($A=0$) with $\rho_c$ chosen high enough to cluster slightly closer to the corners than the Newman poles. Then, we will increase our number of collocation points $N_c$ until the middle of the edges are sufficiently resolved. We note that $\max_{\beta \in (0,2)} \sigma (\beta)=\sqrt{2}\pi$, so a satisfactory choice is $\rho_c = \sqrt{2(m+1)}\pi$. This implies our minimum distance from the corner to a collocation point occurs at $t=1$ with distance $\frac12h_{\min}e^{-\sqrt{2(m+1)}\pi}$, which will generally be closer than the distance to the closest pole $\frac{\sqrt{2}}{2}h_{\min}e^{-\sigma (\sqrt{m+1}-1)}$ for simple enough geometries.}

\subsection{Numerical inversion of the Laplace transform.}\label{sec:LT}

\begin{figure}[htbp]\label{fig:talbotcontour}
    \centering
    \subfigure[Talbot contours for $t=0.1$ and varied $M$.]{\includegraphics[width=0.475\textwidth]{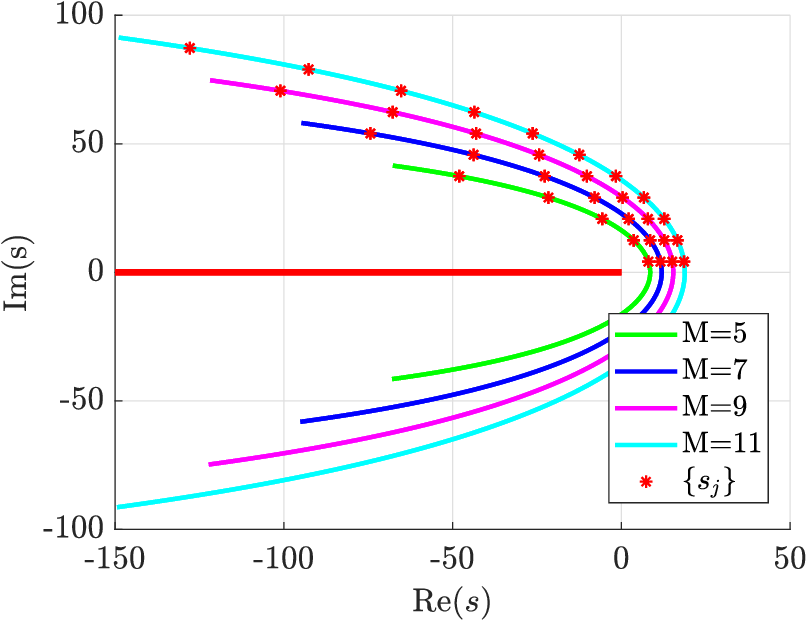}} \qquad  
    \subfigure[Talbot contours for $M=9$ and varied $t$.]{\includegraphics[width=0.475\textwidth]{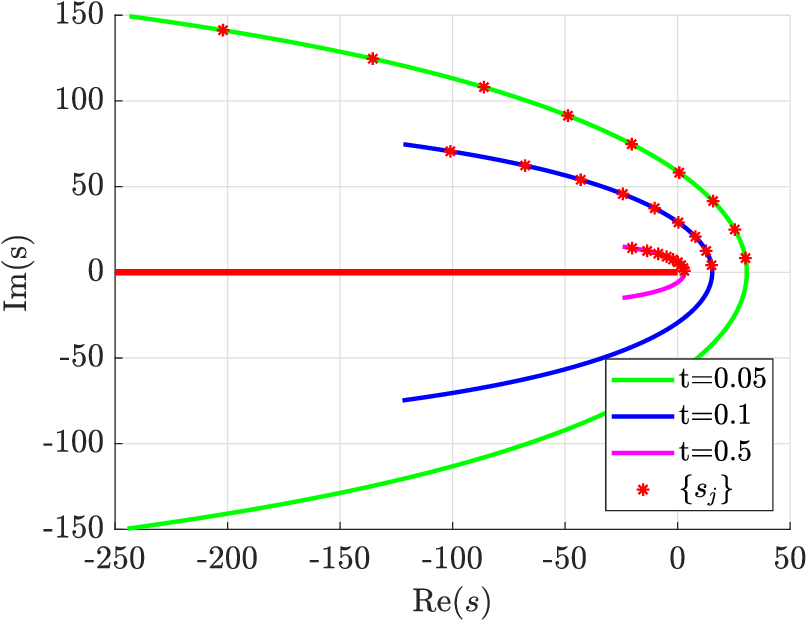}}
    \caption{The Talbot contour \eqref{eq:Talbot} and the mid-point quadrature points $\{ s_j\}$. Shown are contours and their discretizations for fixed $t=0.1$ and varying $M$ (left) and varying $t$ with fixed $M=9$ (right). \label{fig:talbot}}
\end{figure}

In this section we describe our process for inverting the Laplace transform to return to time \al{domain}. The inverse Laplace transform is defined by the Bromwich integral, given by
\begin{equation}\label{eq:invLT}
u(z;t) = \mathcal{L}^{-1}(\hat{u}) = \frac{1}{2 \pi i} \int_B \hat{u}(z;s)e^{st} \diff s,
\end{equation}
where $B = \{ \gamma + i s | s \in \mathbb{R}\}$ is the Bromwich contour. The value $\gamma\in\mathbb{R}$ is selected so that all poles of $\hat{u}(z;s)$ with respect to $s$ are to the left of $\mbox{Re}(s)= \gamma$. The  integrand of \eqref{eq:invLT} tends to be oscillatory along $B$ which can require a lot of terms to compute. However, since our problem has poles that lie solely along the negative real axis, we can deform the contour of integration into the left hand plane so that the integrand of \eqref{eq:invLT} decays rapidly. The optimal choice of contour should balance moving into the left half-plane rapidly while not getting too close to the singularities, nor too far such that the exponential factor in \eqref{eq:invLT} becomes very large. A highly optimized choice, known as the modified Talbot contour \cite{Weideman2015}, is given by the general form
\bsub\label{eq:Talbot}
\begin{equation}\label{eq:Talbot_a}
B_T=  \Big\{ \frac{2M}{t} \rho(\theta)  \, | -\pi < \theta < \pi \, \Big\}, \qquad \rho(\theta) = \big(-\sigma + \mu\theta \cot(\alpha\theta) + \nu i \theta\big),
\end{equation}
 where $2M$ is the number of quadrature points applied to approximate the integral (see Fig.~\ref{fig:talbot}). A saddle point analysis \cite{Weideman2015} yields the following values which optimize the convergence rate:
 \begin{equation}
     \sigma = 0.6122, \qquad \mu = 0.5017, \qquad \alpha = 0.6407, \qquad \nu = 0.2645.
 \end{equation}
\esub
 The parametrization of $B_T$ given in \eqref{eq:Talbot}, expresses the integral as
\begin{equation}\label{eq:Talbot_b}
    u(z;t) = \frac{M}{ \pi i t} \int_{-\pi}^\pi \hat{u}\Big(z;\frac{2M\rho(\theta)}{t}\Big) e^{2M\rho(\theta)} \rho'(\theta) \diff \theta = \text{Re}\left[\frac{2M}{ \pi i t} \int_{0}^\pi \hat{u}\Big(z;\frac{2M\rho(\theta)}{t}\Big) e^{2M\rho(\theta)} \rho'(\theta) \diff \theta\right].
\end{equation}
In the second step of \eqref{eq:Talbot_b}, the symmetry of \eqref{eq:Talbot} has been used to express the integral as twice the real part of the value along the upper branch. The contour \eqref{eq:Talbot} is particularly effective for the present case of parabolic problems whose singularities lie along the negative real axis. The theoretical convergence rate is $\mathcal{O}(10^{-1.2M})$. Our recovery of this convergence rate is presented in Fig.~\ref{fig:LTconvergence}.

Given a fixed $M$, our approximation method is as follows: Discretize $(0,\pi)$ at equally spaced points $\theta_1 ,\ldots,\theta_M$ with step-size $h = \pi/M$. For a given $t>0$, define $\omega (\theta) = e^{2M\rho (\theta)} \rho'(\theta)$ and $s(\theta) = 2M\rho(\theta)/t$. We then evaluate weights $\omega_j = \omega(\theta_j)$ and points $s_j = s(\theta_j)$ to define the midpoint quadrature rule
\begin{equation}
    u(z;t) =  \text{Re}\Big(\frac{2}{it} \sum_{j=1}^M \omega_j \hat{u}(z;s_j) \Big), \qquad \theta_j = \Big(j-\frac{1}{2}\Big)\frac{\pi}{M},\qquad j = 1,\ldots,M.
\end{equation}
We remark that this choice of contour has the advantage of converging very rapidly\al{,} with errors reduced to near machine precision using $M=9$ transform evaluations. Additionally, as the numerical examples of Sec.~\ref{sec:Results} will demonstrate, this integration method is robust over a wide range of times $t\in(0,\infty)$. A disadvantage of this approach is that it incurs a significant computational cost to recompute the transforms $\hat{u}(z;s_j)$ for each specified value of $t$. In the conclusion, we discuss methods for ameliorating this burden.

%
%
%
%
%

\section{Comparable methods for numerical validation and application.}\label{sec:compmethods}
We need various ways to validate and benchmark our solution. On most domains, analytical solutions \al{are unavailable}, and hence we will validate with other numerical and asymptotic methods. \hl{We will validate with a high-order boundary integral method, to ensure our solutions are accurate beyond low boundary error. In addition, we will compare with two lower-accuracy methods, a particle based Kinetic Monte Carlo (KMC) method and an asymptotic method based on the limit of small, well separated absorbing bodies \cite{che-lin-her-qua2022}. The surface flux $j(t) = \int_{\partial \Omega} D\partial_n u \diff s$ of our approximation has an integrand that can be acquired analytically, that is then numerically integrated and compared with out other methods.}

\subsection{\hl{Boundary Integral Equation for the Modified Helmholtz equation}}

\hl{In the boundary integral equation, we represent our solution using the density of a layer potential, in our case the double layer potential
\begin{equation}
	\mathcal{D}[\sigma_D](z)=\frac1{2\pi} \int_{\partial \Omega} \partial_{n}G(z,\xi) \sigma_D(\xi) \diff \xi,
\end{equation}
Where in this case $n$ is the inward facing normal. The integral equation corresponding to the exterior Dirichlet problem with initial condition $u_0(z)=\delta(z-z_0)$ is then \cite{colton2013integral}
\begin{equation}
	\frac12 \sigma_D(z)-\mathcal{D}[\sigma_D] (z) = -G(z,z_0). \label{eq:mhebie}
\end{equation}
We will use the integral equation package \texttt{chunkIE} \cite{askham2024chunkie}, which has a built in solver for the density of a double layer potential for the standard Helmholtz equation with kernel $G(z,\xi)=\frac{i}4H^{(1)}_0(k|z-\xi|)$, where $H^{(1)}_0(z)$ is the Hankel function of the first kind. The identity \cite{DLMF}
\begin{equation}
	K_0(z)=\frac{\pi i}2 H^{(1)}_0(iz),
\end{equation}
is used to represent the solution of \eqref{eq:Helmholtz} with the substitution of $k=i\sqrt{s}$.
}

\subsection{Kinetic Monte Carlo (KMC) methods} In the particular case of solving \eqref{eq:main} when $f\equiv0$ and the initial condition is given by a Dirac source, $u_0(z) = \delta(z-z_0)$, we corroborate with a complementary method for solving the PDE \eqref{eq:main} through sampling of the Langevin equation \eqref{eq:prob_dists}. The KMC method is a technique that massively accelerates the sampling of \eqref{eq:prob_dists} by breaking down random paths into simpler geometric steps, or projectors, that can be solved exactly and tabulated for rapid access. By carefully combining several projection steps, it is possible to create exact simulations in arbitrary geometries in two and three dimensions \cite{OKMascagni2004,che-lin-her-qua2022,lindsay2024,LBS2018,ye2024}. The simulation of the paths \eqref{eq:prob_dists_a}, and the statistics of arrivals of particles to $\partial\Omega$, gives Monte-Carlo estimates for the quantities \eqref{eq:prob_dists_b} which we can compute from our LM solution of \eqref{eq:main}. This provides a highly non-trivial method for validating our results over wide timescales.

\subsection{Matched asymptotic expansion solution for multiple bodies.}\label{sec:MAA}
The LM method can trivially be extended to multiple bodies by adding collocation points together with corresponding Newman and Runge expansions to the series approximation. This allows us to approximate the solution of \eqref{eq:main} around multiple bodies. A useful validation of this approximation will be a recent matched asymptotic solution of \eqref{eq:main}, derived in the limit of small well-separated bodies \cite{che-lin-her-qua2022} which we briefly describe. \al{The $N_B$ bodies $\{\Omega_{k}\}_{k=1}^{N_B}$ are centered at $\{z_k\}_{k=1}^{N_B}$ with a common scale factor $\varepsilon$, such that
\[
\Omega_{k} = z_k + \varepsilon \, \mathcal{A}_k, \qquad k = 1,\ldots,N_B.
\]
Here $\mathcal{A}_k$ describes the shape of the $k^{th}$ body.} It was recently shown \cite{che-lin-her-qua2022} that the solution of \eqref{eq:Helmholtz} is of form
\bsub\label{eqn:asy_solution}
\begin{equation}\label{eqn:asy_solution_a}
    \hat{u}(z;s) = G(z,z_0;s) - \sum_{k=1}^{N_B} \hat{j}_k (s) G(z,z_k;s) + \mathcal{O}(\eps), \qquad \nu_k = \frac{-1}{\log \eps d_k},
\end{equation}
\al{in the limit as $\eps\to0$.} Here $G$ is the Green's function \eqref{eq:Greens} and $d_k$ is a constant called the \emph{logarithmic capacitance} or \emph{logarithmic capacity} \cite{Venu2015,Polya1951,Baddoo2021} which is determined uniquely by the shape of $\mathcal{A}_k$. \al{The value of $d_k$ is known exactly for some regular geometries (e.g. ellipses, squares and triangles) and can be numerically approximated for general shapes (details in Appendix \ref{sec:log_cap}). The fluxes into each body are defined as $j_k(t) = \int_{\partial\Omega_k}D\partial_n u\, \diff s$ with their corresponding Laplace transforms being $\{\hat{j}_k(s)\}_{k=1}^{N_B}$. These satisfy the linear system}
\al{\begin{equation}\label{sysmain}
 (\mathcal{I} + 2\pi D\, \mathcal{V} \mathcal{G} )\, \hat{\textbf{J}} = 2\pi D \mathcal{V} \textbf{g}_0, \qquad [\mathcal{G}]_{i,j}= \left\{ \begin{array}{rl} R(s),\quad i= j, \\[5pt]  G(z_i,z_j;s),\quad i\neq j. \end{array} \right. \qquad [\mathcal{V}]_{i,j}= \left\{ \begin{array}{rl} \nu_i \quad i= j, \\[5pt]  0,\quad i\neq j, \end{array} \right. 
\end{equation}
where $\mathcal{I}\in\mathbb{R}^{N_B\times N_B}$ is the identity matrix,
$\hat{\textbf{J}}\in\mathbb{R}^{N_B}$ and $ \textbf{g}_0\in\mathbb{R}^{N_B}$ are given by
\begin{align}
\hat{\textbf{J}}= [\hat{j}_1,\hat{j}_2,\ldots,\hat{j}_{N_B}]^T,\qquad
 \textbf{g}_0=[G(z_1,z_0;s), G(z_2,z_0;s), \ldots,G(z_{N_B},z_0;s)]^T.
\end{align}}
\esub
The matrix $\mathcal{G}$ describes the spatial interactions between the bodies while the vector $\textbf{g}_0$ describes the
influence of the initial location on each of the bodies. Here $R(s)$ is the regular part of $G(z,\xi;s)$ at the source. The small argument asymptotics $K_0(z)\sim-\log(z)+\log2-\gamma_e$ as $z\to0$ give this term as
\begin{equation}
 R(s)= \lim_{z\to \xi} \Big(G(z,\xi;s) + \frac{1}{2\pi D} \log|z-\xi| \Big) = \frac{1}{2\pi D}\Big(\log2-\gamma_e-\log\sqrt{s/D}\, \Big),
\end{equation}
where $\gamma_e\approx 0.5772$ is the Euler-Mascheroni constant. Hence, given the centering coordinate $z_k$ and capacitance $d_k$ of each body, we solve the linear system \eqref{sysmain} and apply the Talbot method (cf.~Sec.~\ref{sec:LT}) to obtain the time-dependent fluxes $j_k(t)$ and corresponding cumulative fluxes $c_k(t) = \int_0^tj_k(\tau) \diff \tau$.

\section{Numerical Results.}\label{sec:Results}
In this section, we provide a variety of examples that validate our method for single and multiple bodies. \hl{As our solutions are built from solutions to \eqref{eq:Helmholtz}, we will validate it first. Then,} we will \hl{examine} the effect of varying the number of transform evaluations $M$ in the inverse Laplace transform procedure, eventually adopting $M=9$ as the default parameter. In each of the examples shown, the diffusivity takes the value $D=1$. We will examine the accuracy of the method by determining the errors
\bsub
\begin{align}
    \label{eq:rel_errors_a} \mathcal{E}^{\infty}[\hat{u}] = \sup_{z\in\partial\Omega}\|\hat{u}(z;s) - f(z)/s\|,\\[4pt]
    \label{eq:rel_errors_b}\mathcal{E}^{\infty}[u] = \sup_{z\in\partial\Omega}\|u(z,t) - f(z)\|.
\end{align}
\esub
 Recall that $\hat{u}-f(z)/s=\hat{u}_h+\hat{u}_p-f(z)/s=\hat{u}_h-(f(z)/s-\hat{u}_p)=\hat{u}-\tilde{f}$, so the error is equivalent to the error between $\hat{u}_h$ and $\tilde{f}$ from \eqref{eq:uh_c}. Accordingly, when taking \al{the} relative error, we scale by the factor $\sup_{z \in \partial \Omega} \|\tilde{f}(z;s)\|$ in the Helmholtz problem, and the equivalent in the heat problem. \hl{We approximate these errors by oversampling on $\partial \Omega$ with distribution $e^{-\left( \sqrt{2(m+1)}\pi+\log(10) \right)t}$ from a corner to a midpoint, using $3$ times as many sample points than collocation points.}
 
 Unless stated otherwise, our primary demonstrations will set $f(z)\equiv0$, and $u_0(z)=\delta(z-z_0)$ for some $z_0$, \hl{which solves} the underlying problem for calculating FPT distributions \eqref{eqn:bdflux}. \hl{Note that for this choice of $u_0(z),$ the boundary term \eqref{eq:uh_c} reduces to}
\[
\tilde{f}(z;s) = \frac{f(z)}{s} - \int_{\mathbb{C}\setminus\Omega} u_0(\xi)G(z,\xi;s)\diff \xi = \frac{f(z)}{s} - G(z,z_0;s).
\]
\hl{which is readily evaluated.}  Recalling that our number of degrees of freedom is $2N$, we form our overdetermined linear system by setting the number of collocation points at \hl{$N_c=5\times2N$ (One could get away with a lower scaling factor, but we amplify for robustness)}. Finally, we choose \hl{$N_2=3.5\sqrt{m}$} as the order of our Runge expansions.

\subsection{ \hl{Validation of the Modified Helmholtz solution}}
\begin{figure}[htbp] 
	\centering
	\subfigure[Relative residual.]{\includegraphics[width=0.45\textwidth]{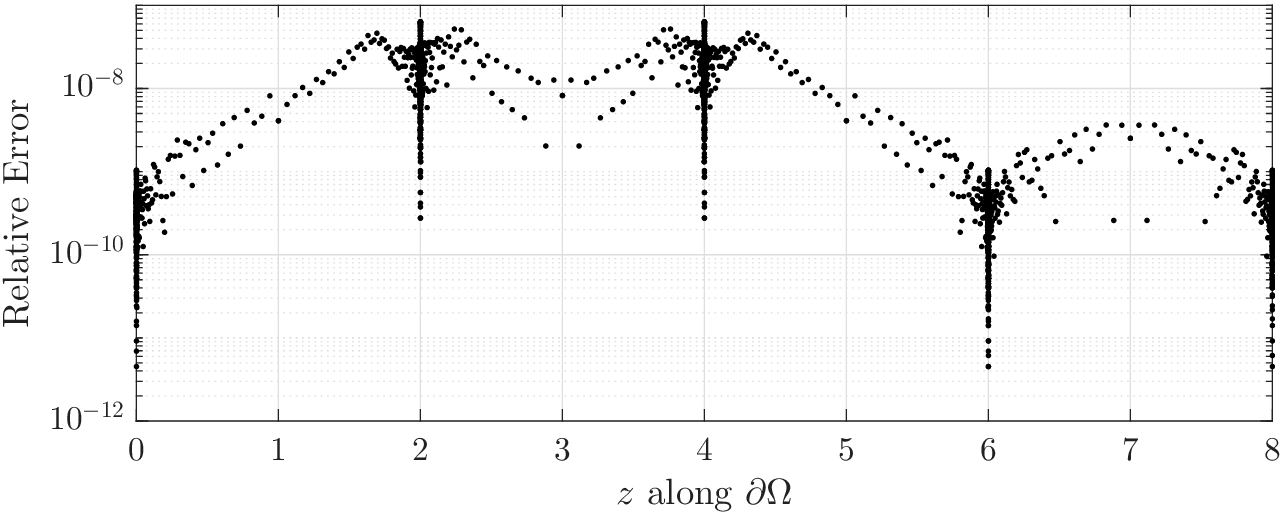}} \qquad
	\subfigure[Relative error on over-sampled grid.]{\includegraphics[width=0.45\textwidth]{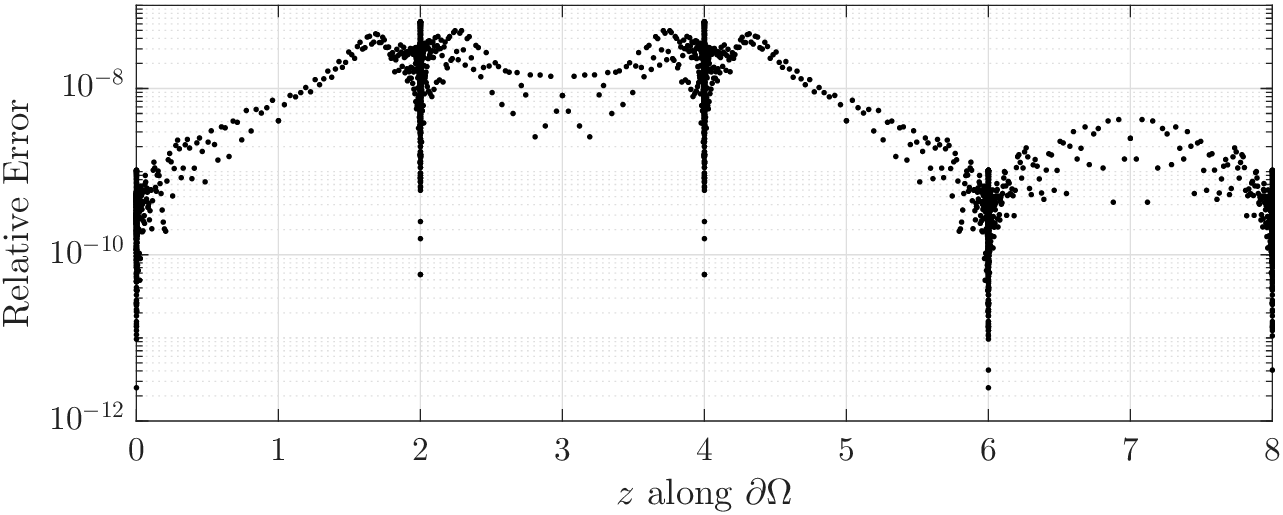}} \\
	\subfigure[Relative residual.]{\includegraphics[width=0.45\textwidth]{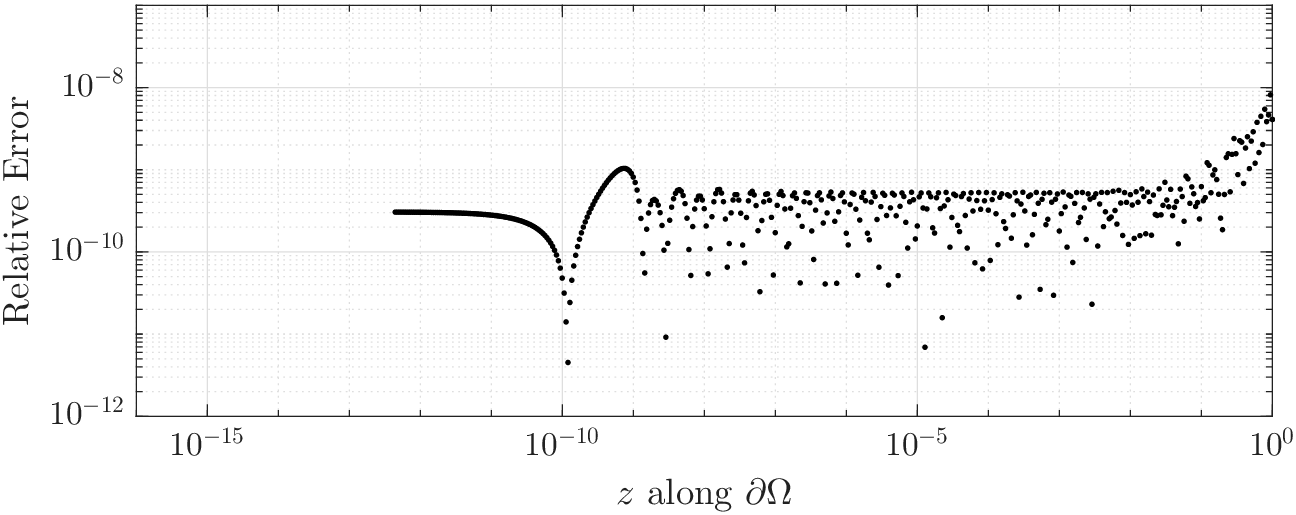}} \qquad
	\subfigure[Relative error on over-sampled grid.]{\includegraphics[width=0.45\textwidth]{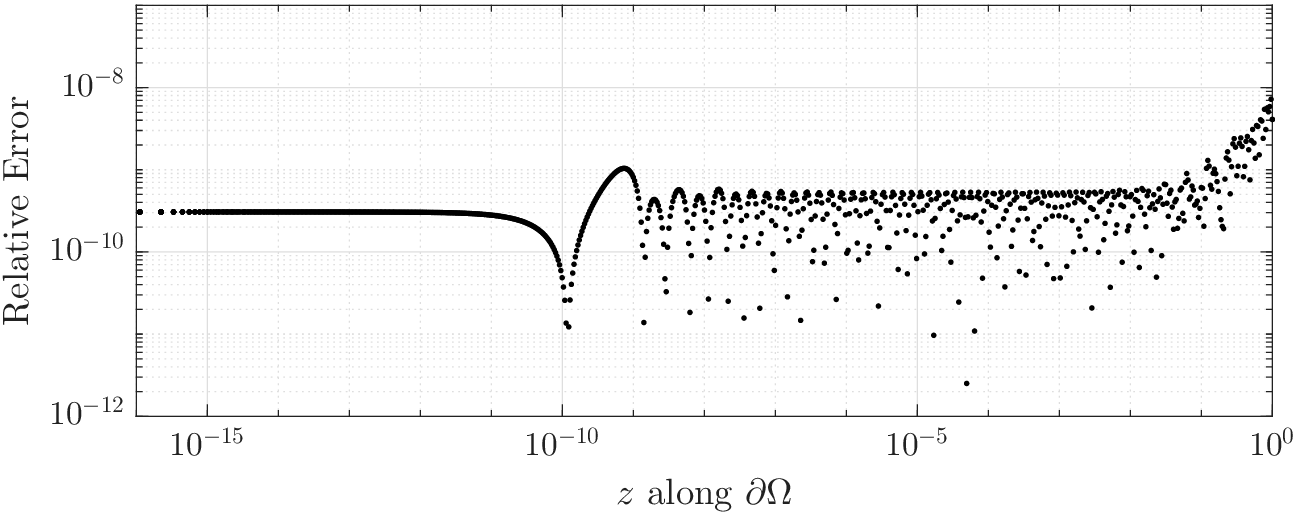}} \\
	\caption{\hl{Point-wise relative residual and over-sampled relative error for $s=-30.16+8.31$, $m=40$ on panels (a-b): a parametrization of the boundary $\partial \Omega$ starting from $-1-1i$ and going counter-clockwise, and panels (c-d): log-scale on $[0,1]$ of the parametrization of the boundary, showing the sample clustering going deeper than the residual to the corner.}} \label{fig:errstruct}
\end{figure}

\hl{We first validate our solution to \eqref{eq:Helmholtz} by plotting the boundary error on an oversampled grid for a square geometry with vertices at $\pm 1 \pm i$ and a point source at $z_0=2$. In Fig.~\ref{fig:errstruct} we plot the relative residual and the relative error on the oversampled grid that $\mathcal{E}^\infty [\hat{u}]$ is calculated from for a lower value of $m$ to show detail. Excellent accuracy is observed all the way to the corner as shown in the zoomed view of Fig.~\ref{fig:errstruct}(d).}

\begin{figure}[h] 
    \centering
    \subfigure[$|\mathrm{Re}(\hat{u}(z,s))|$ at $s=30.16+8.31i$.]{\includegraphics[width=0.31\textwidth]{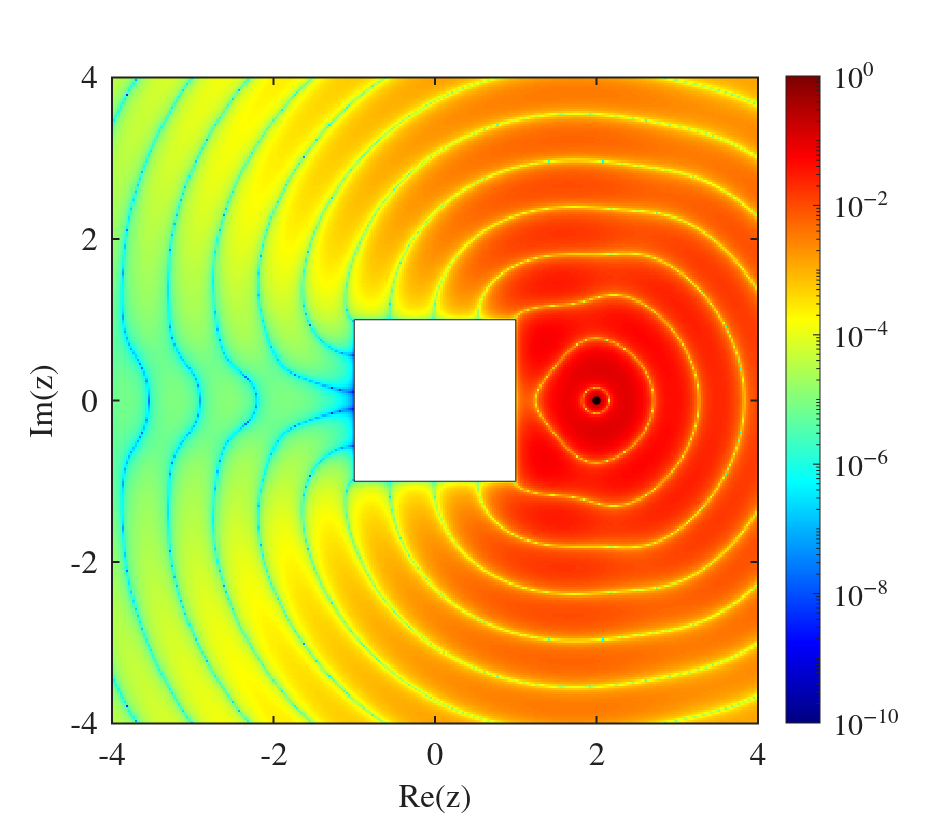}} \quad
    \subfigure[$|\mathrm{Re}(\hat{u}(z,s))|$ at $s=-0.72+58i$.]{\includegraphics[width=0.31\textwidth]{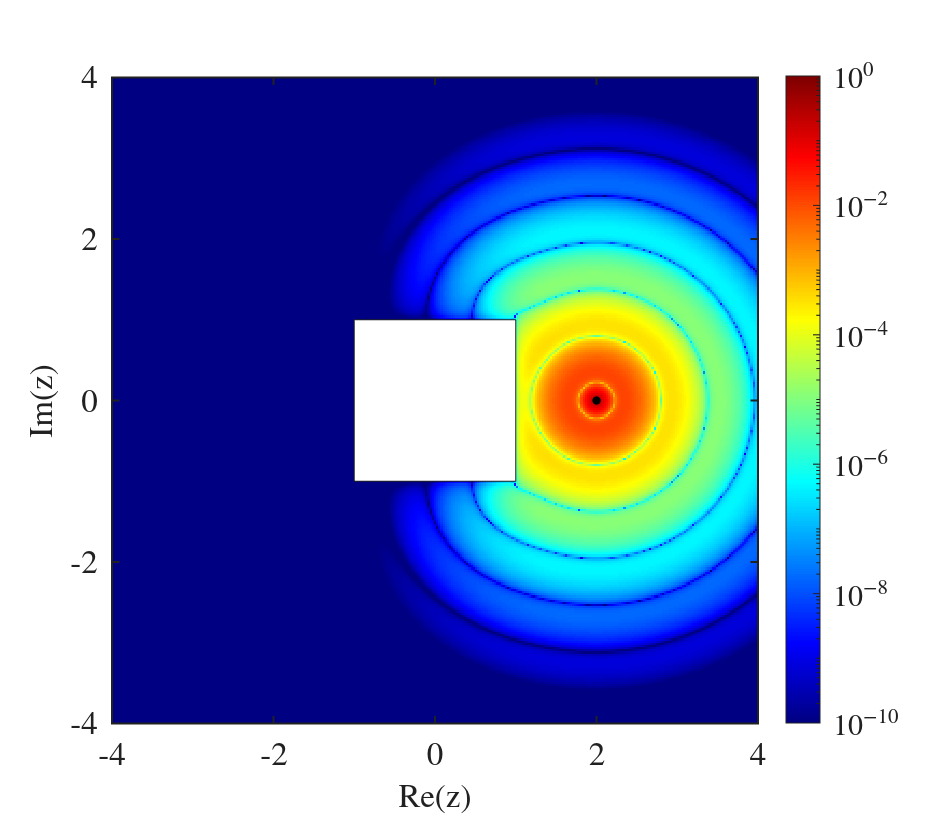}} \quad  
    \subfigure[$|\mathrm{Re}(\hat{u}(z,s))|$ at $s=-86.09+108i$.]{\includegraphics[width=0.31\textwidth]{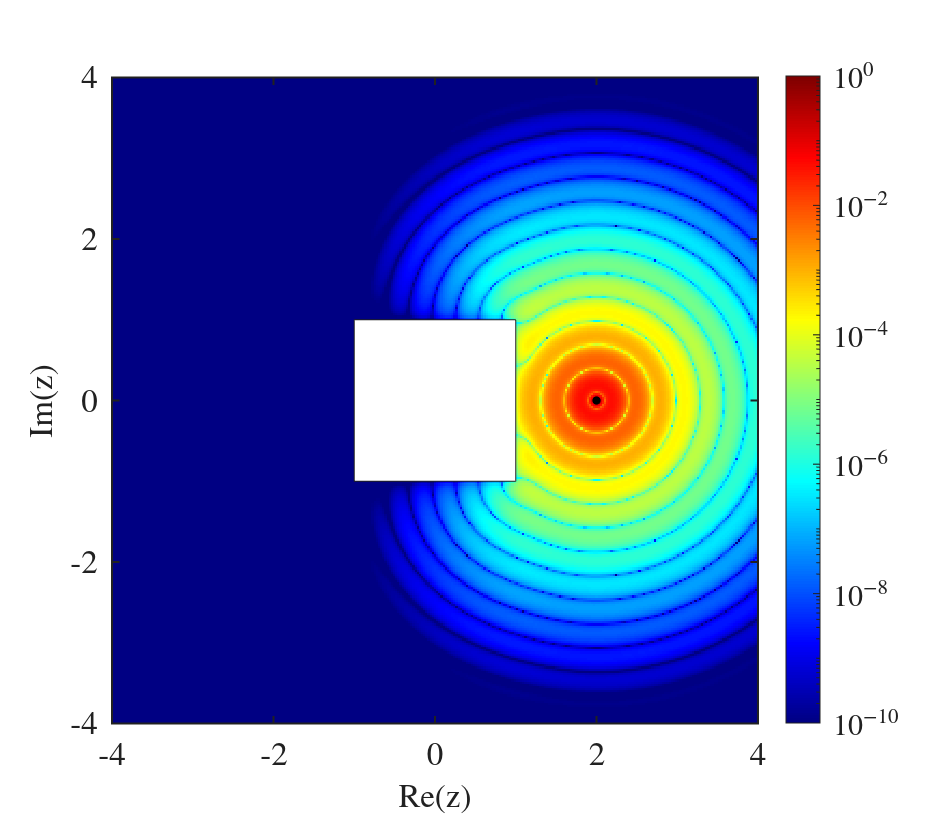}}
    \subfigure[Convergence of errors for $\hat{u}(z,s)$.]{\includegraphics[width=0.475\textwidth]{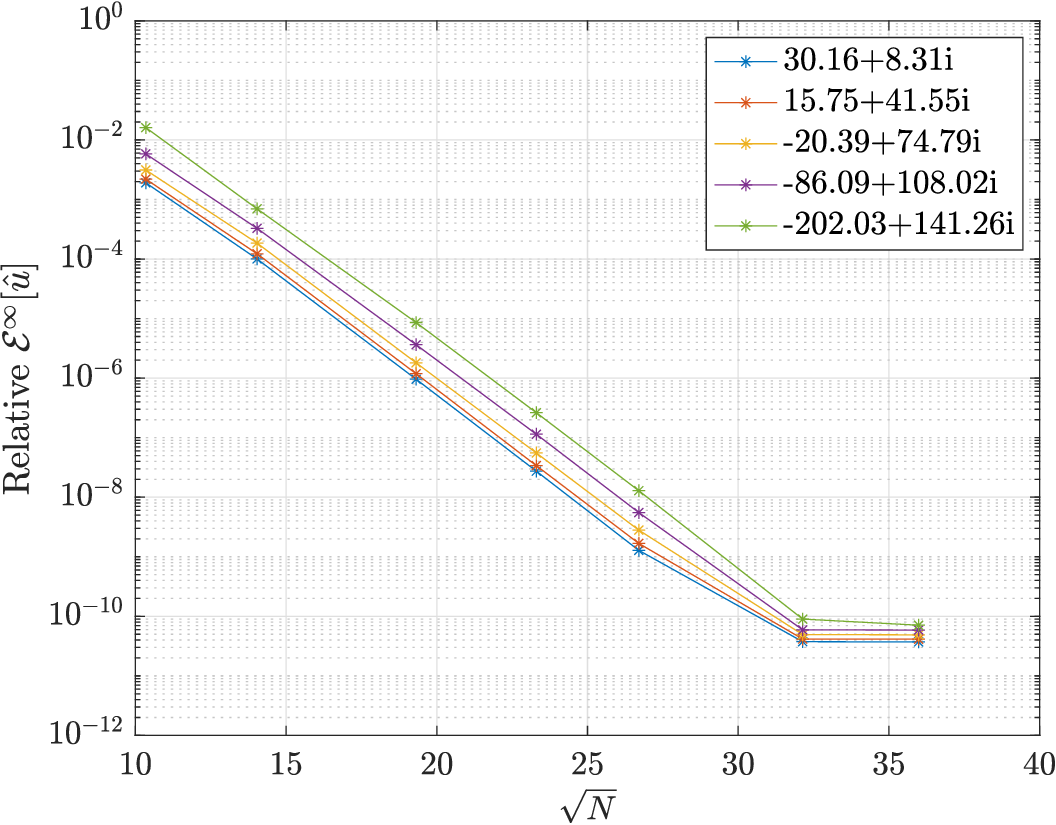}} \qquad
	\caption{Panels (a-c): Approximate solution of \eqref{eq:Helmholtz} for various $s$ in the Talbot contour of $t=0.1$, with \hl{$m=90$}. $\mathcal{E}^\infty[\hat{u}]$ for each given by (a) \hl{ $2.02\times 10^{-10}$, (b) $8.24\times 10^{-11}$, (c) $3.58\times 10^{-10}$.} (d): Convergence of the relative errors for various $s$ in the Talbot contour for time $t=0.1$ in the same configuration as (a-c).}\label{fig:squaremhh}
\end{figure}

In Fig.~\ref{fig:squaremhh}, we plot solutions of \eqref{eq:Helmholtz} for several values of $s$ corresponding to the Talbot contour of $t=0.1$. \hl{In Fig.~\ref{fig:squaremhh}, we display results of a convergence study for the square example by calculating $\mathcal{E}^{\infty}[\hat{u}]$ from \eqref{eq:rel_errors_a} for an increasing number of terms $N$. We observe the convergence rate to be linear against $\sqrt{N}$ which demonstrates root-exponential convergence, and an attainment of relative error $\approx 10^{-10}$ for a range of $s$ values associated with the inverse Laplace transform process. Eventually the convergence stagnates, which has been observed in the LM literature and is due to the process of solving with an overdetermined basis \cite{Gopal19B}.
 }

\begin{figure}[h] 
    \centering
    \subfigure[$|\mathrm{Re}(\hat{u}(z,s))|$ via Boundary Integral Method]{\includegraphics[width=0.45\textwidth]{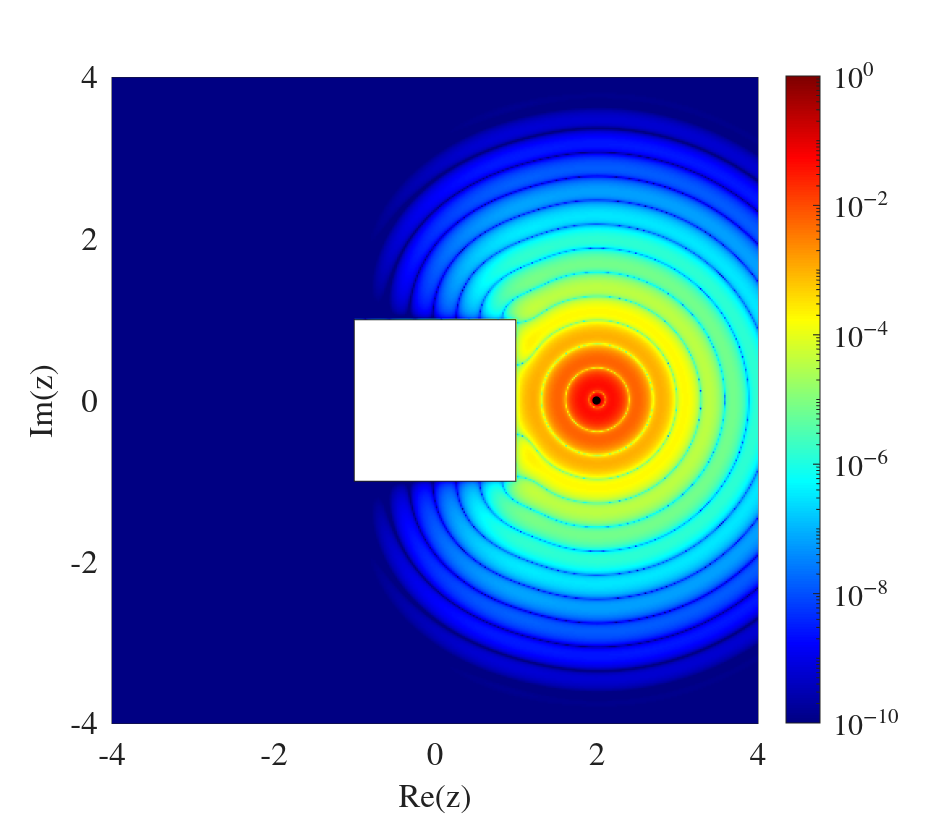}} \quad  
    \subfigure[Absolute difference in BIE and LM solutions]{\includegraphics[width=0.45\textwidth]{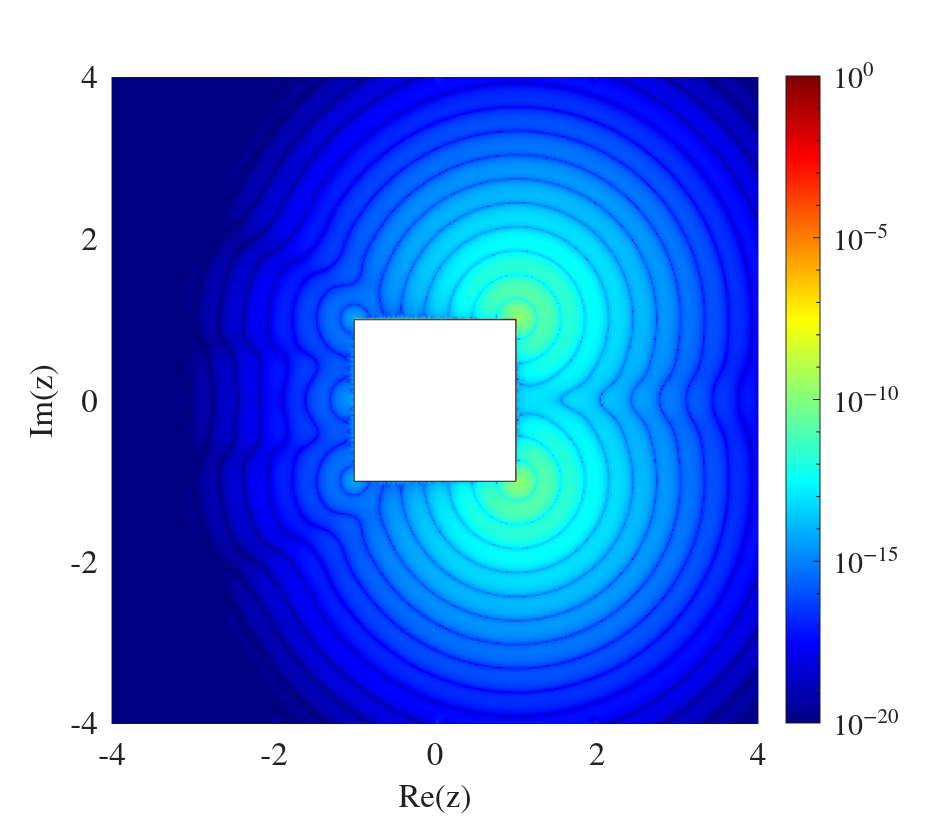}}
	\caption{\hl{Comparison with a boundary integral method with the configuration in Fig.~\ref{fig:squaremhh} for single $s=-86.09+108i$. (a) Solution of \eqref{eq:Helmholtz} via BIE \eqref{eq:mhebie} using \texttt{chunkIE}. A \texttt{chunkgraph} is created with the squares vertices and refined until 64 chunks are on each side.}}\label{fig:BIEcomp}
\end{figure}

\hl{In Fig.~\ref{fig:BIEcomp}, we replicate the result shown in Fig.~\ref{fig:squaremhh}(c) with \texttt{chunkIE} by evaluating the boundary integral equation \eqref{eq:mhebie}. On a logarithmic scale, the solution plots are nearly identical, but when the absolute value of the difference between the two solutions is plotted, we observe a maximum absolute error of $4.47 \times 10^{-10}$ and a maximum relative error around $1.27 \times 10^{-6}$.}

\subsection{Validation of the Talbot integral \& the heat solution}

\begin{figure}[h] 
    \centering
    \subfigure[$M=3$]{\includegraphics[width=0.45\textwidth]{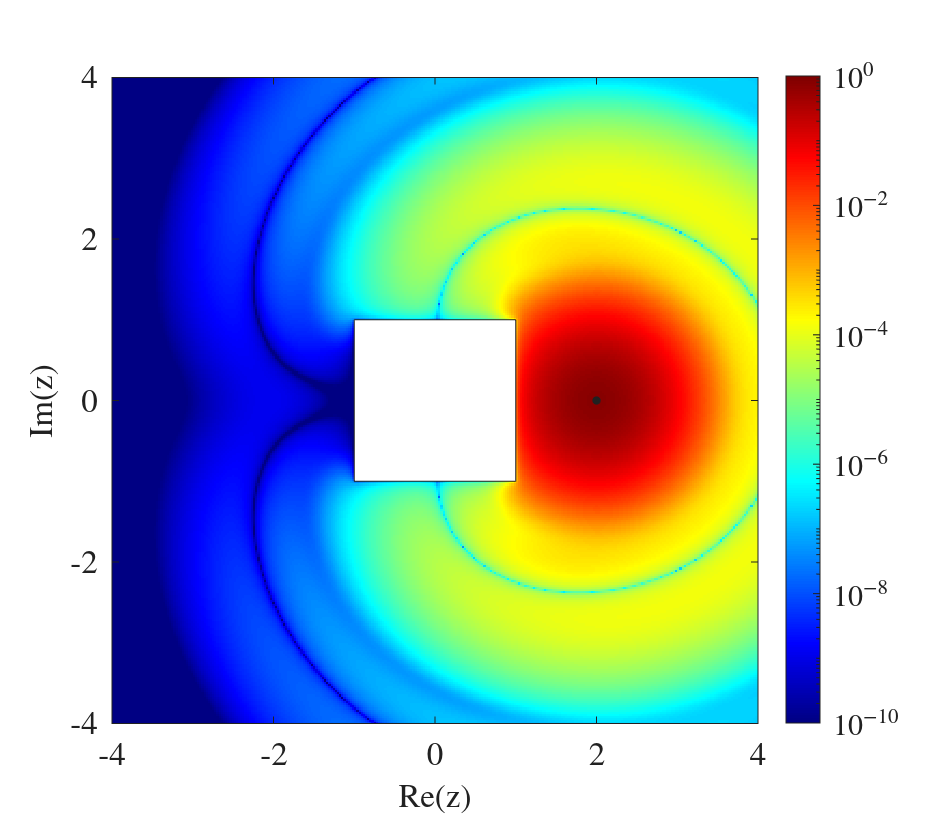}} \qquad
    \subfigure[$M=6$]{\includegraphics[width=0.45\textwidth]{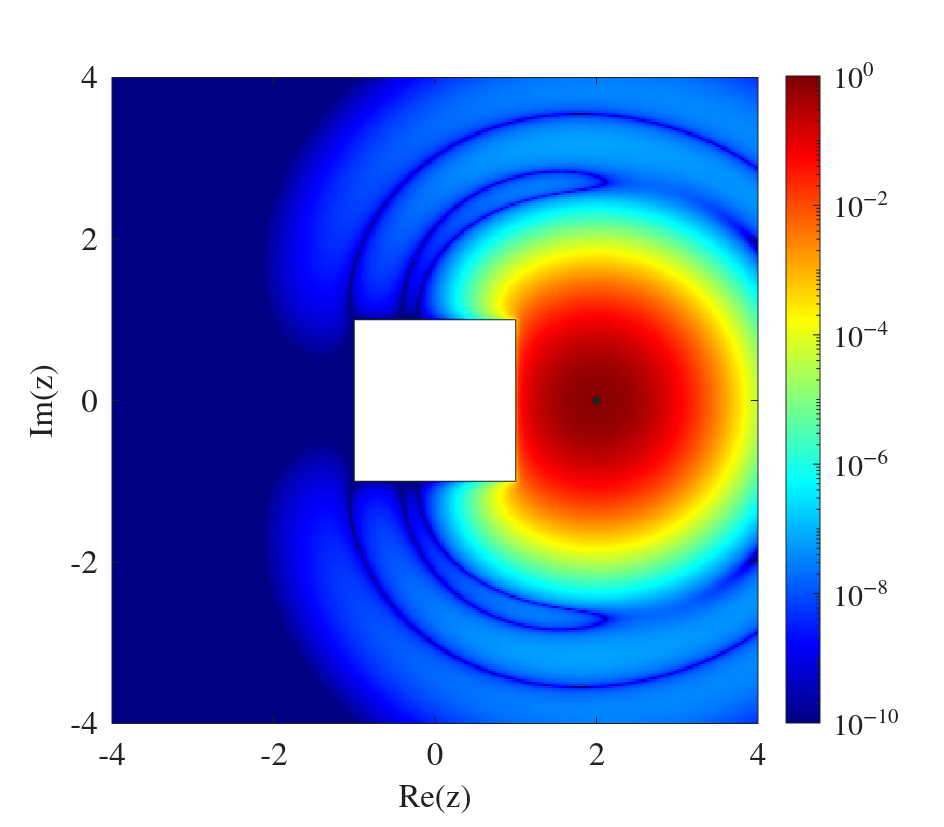}}\\ 
    \subfigure[$M=9$]{\includegraphics[width=0.45\textwidth]{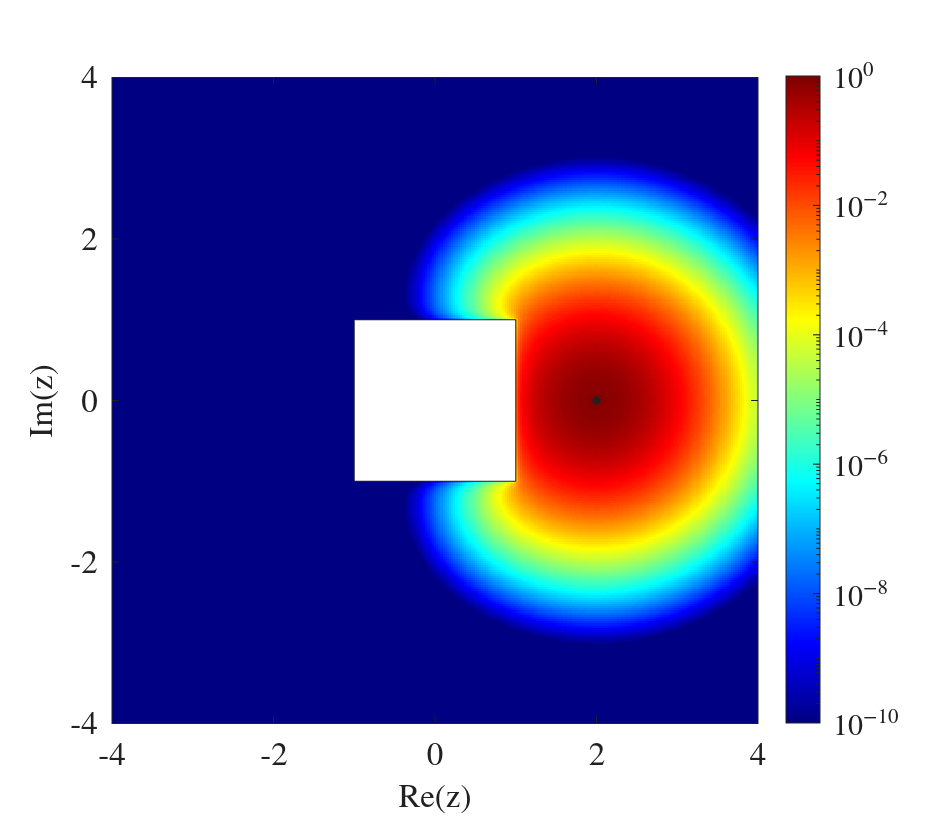}}\qquad
    \subfigure[$M=12$]{\includegraphics[width=0.45\textwidth]{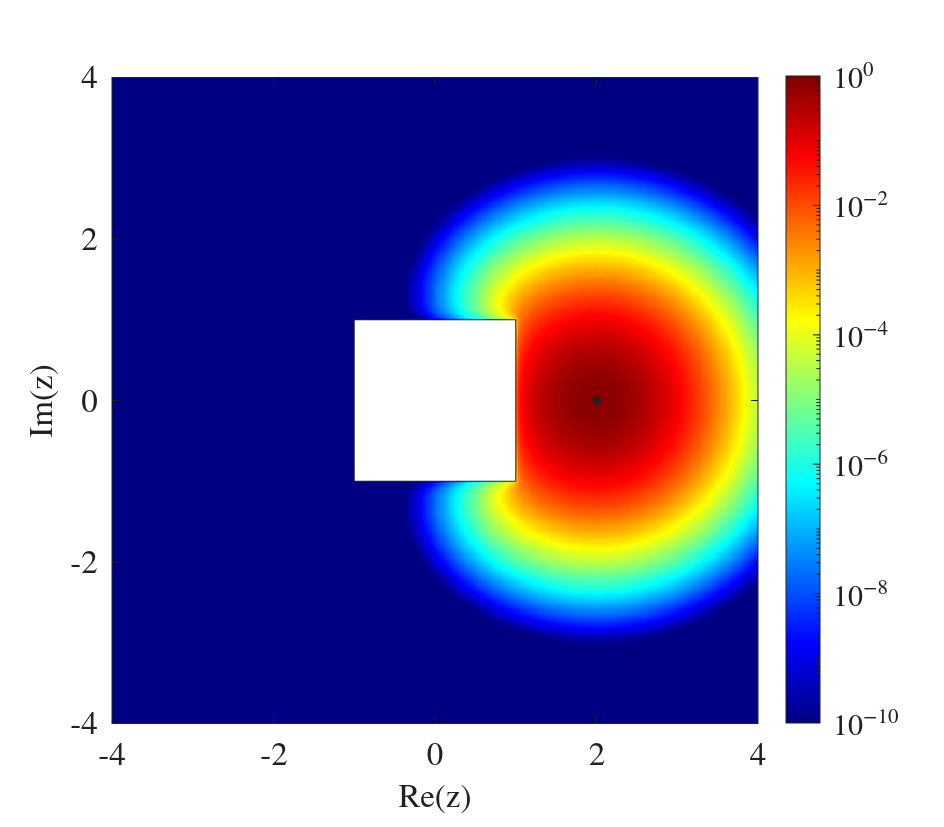}}
    \caption{Solution $|u(z,t)|$ of \eqref{eq:main} at $t=0.1$ with various number $M$ of Laplace transform evaluations. Increasing $M$ removes spurious oscillations in the solution. \label{fig:varyM}}
\end{figure}

In this \hl{section}, we perform validations and examine the convergence properties of the inverse Laplace transform method. Specifically, we seek to demonstrate the appropriate number $M$ of transform evaluations to apply in the Talbot quadrature equation \eqref{eq:Talbot_b}. In Fig.~\ref{fig:varyM}, we display solutions $u(z,t)$ of the system \eqref{eq:main} \hl{for the geometry configuration in Fig.~\ref{fig:squaremhh} at $t=0.1$}. We fix our solution with $m=90$ Newman poles per corner, \hl{matching the Helmholtz solutions and} so that the errors are dominated by the Talbot quadrature. When plotted on a logarithmic scale, the solutions reveal (cf.~Fig.~\ref{fig:varyM}(a,b)) exponentially small spurious oscillations when $M=3,6$ which are resolved (cf.~Fig.~\ref{fig:varyM}(c,d)) for $M=9,12$.

\begin{figure}[h] 
    \centering
	\subfigure[Convergence with respect to $M$.]{\includegraphics[width=0.475\textwidth]{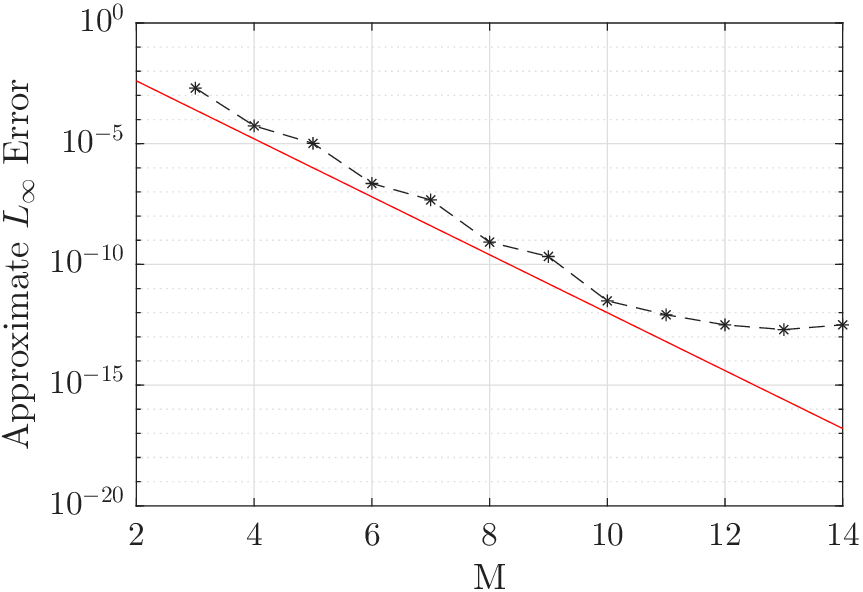}}
	\qquad
    \subfigure[Convergence of errors for $u(z,t)$.]{\includegraphics[width=0.475\textwidth]{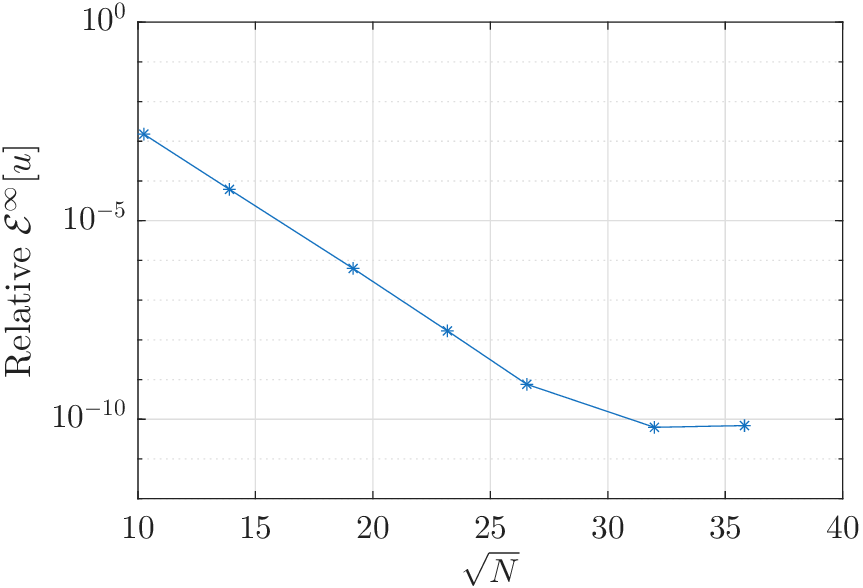}}
	\caption{(a) Convergence of the Talbot method with respect to the number $M$ of transform evaluations. The approximate $L_\infty$ error against $M$ agrees with the theoretical convergence rate \al{of approximately} $\mathcal{O}(10^{-1.2M}$) (straight red line). (b) \hl{Convergence of the relative errors in the time domain for the  solution of equation \eqref{eq:main} at $t=0.1$ in the case of a single square absorber.} \label{fig:LTconvergence} \label{fig:cvgstudiesheat}}
\end{figure}

In Fig.~\ref{fig:LTconvergence}(a), we approximate the $L_\infty$ error by finely sampling our solution in the bulk, and taking the maximum error between our solution for a given $M$ and our solution for $M=15$. We observe the relative error decreasing as a function of transform evaluations $M$, in agreement with the theoretical convergence rate \al{derived in} \cite{Weideman2015}. The choice of $M=9$ is observed to push the error in the bulk close to $10^{-10}$, which we will adopt as our target \hl{(absolute)} accuracy for our other numerical demonstrations. 

\hl{With the value $M=9$ chosen for match the target accuracy of our LM solver, we repeat the convergence study of Fig.~\ref{fig:squaremhh} but for the heat solution at $t=0.1$ . In Fig.~\ref{fig:cvgstudiesheat}, we display the error $\mathcal{E}^{\infty}[u]$ defined in \eqref{eq:rel_errors_b} and again observe root-exponential convergence with a lowest error around $\mathcal{E}^{\infty}[u] \approx 10^{-10}$.} 

\subsection{Example: Simulation for a single triangular absorber.}

\begin{figure}[htbp] \label{fig:varyt}
    \centering
    \subfigure[$u(z,t)$ at $t=0.01$.]{\includegraphics[width=0.31\textwidth]{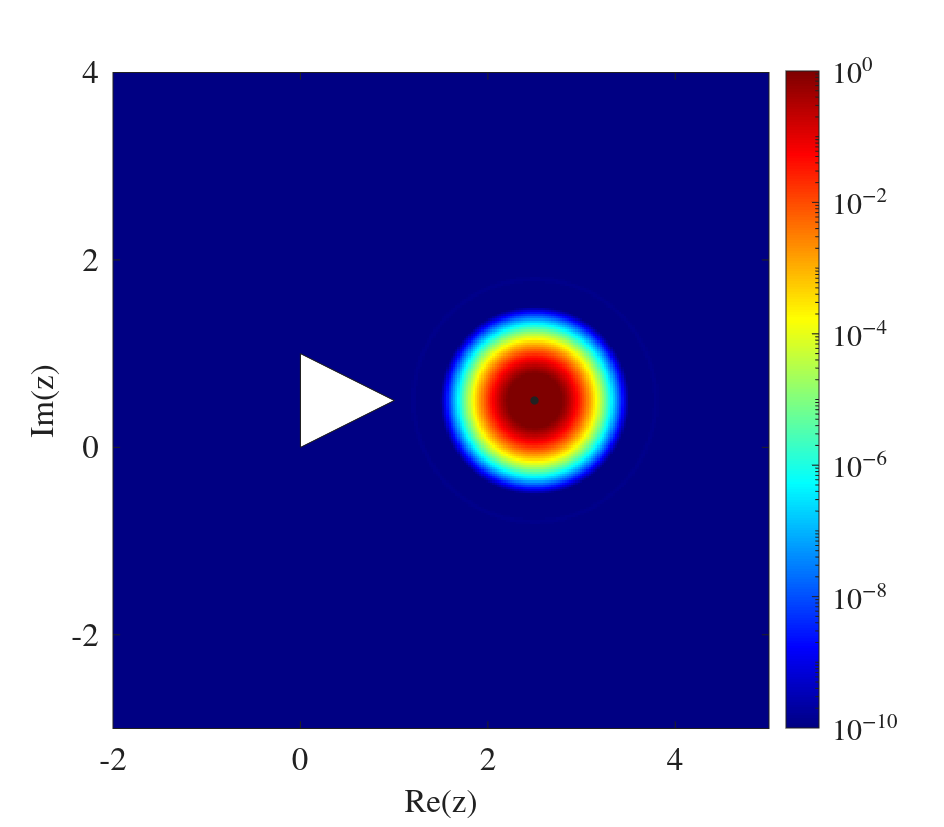}} \quad
    \subfigure[$u(z,t)$ at $t=0.04$.]{\includegraphics[width=0.31\textwidth]{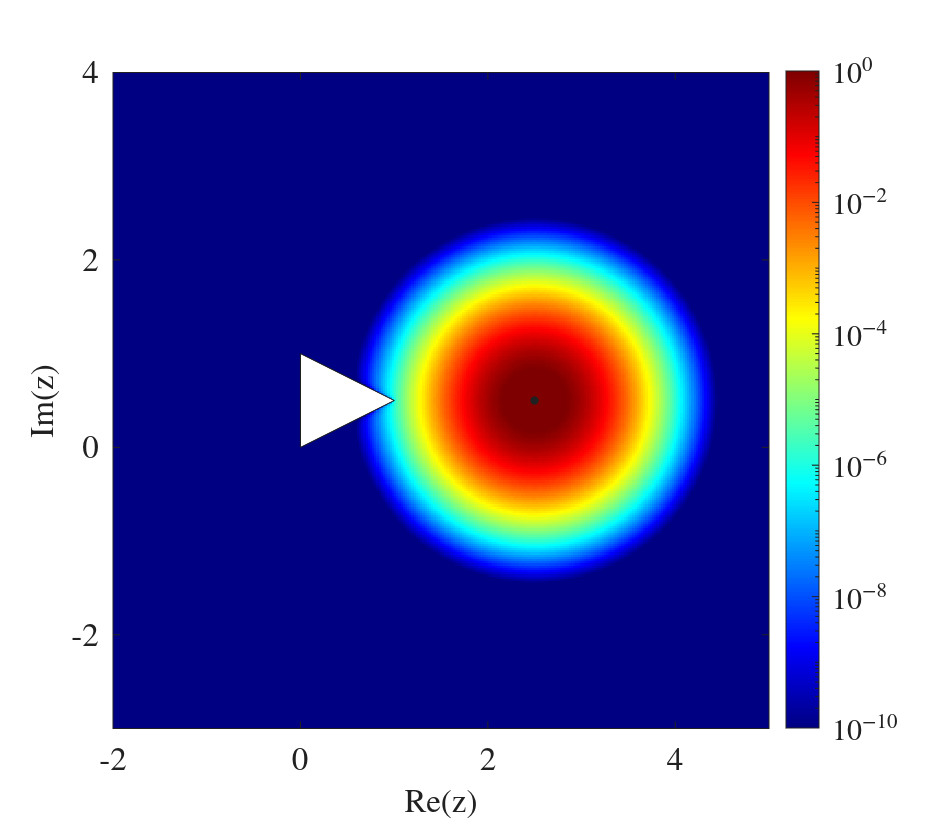}} \quad  
    \subfigure[$u(z,t)$ at $t=0.14$.]{\includegraphics[width=0.31\textwidth]{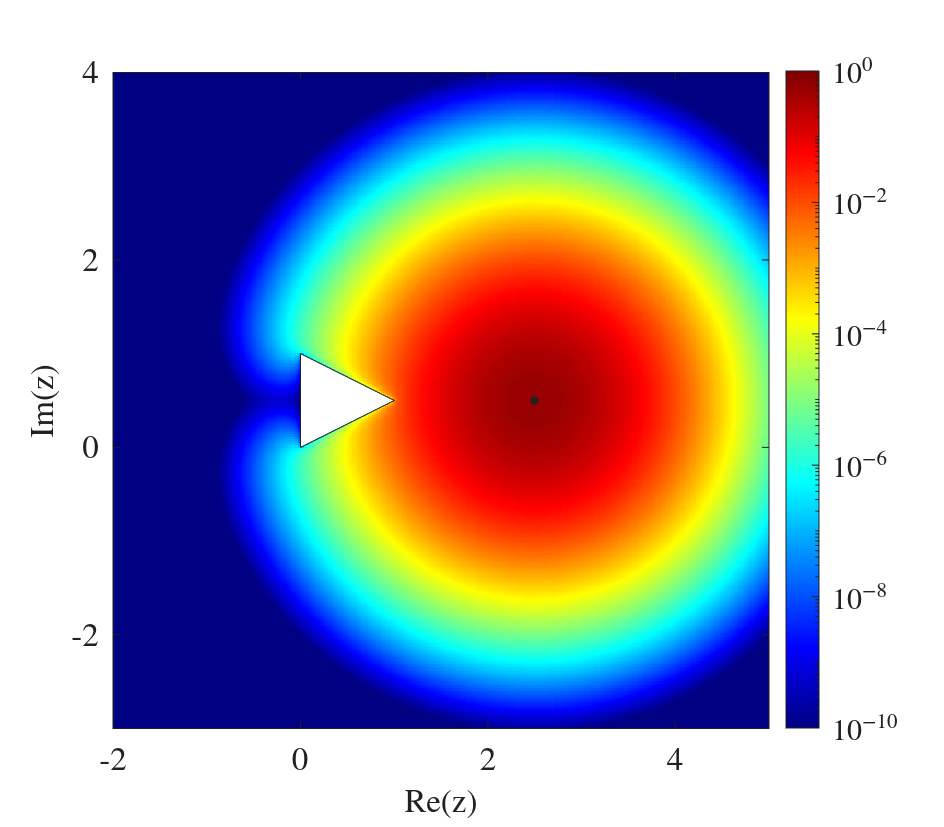}}\\
    \subfigure[$u(z,t)$ at $t= 0.53$.]{\includegraphics[width=0.31\textwidth]{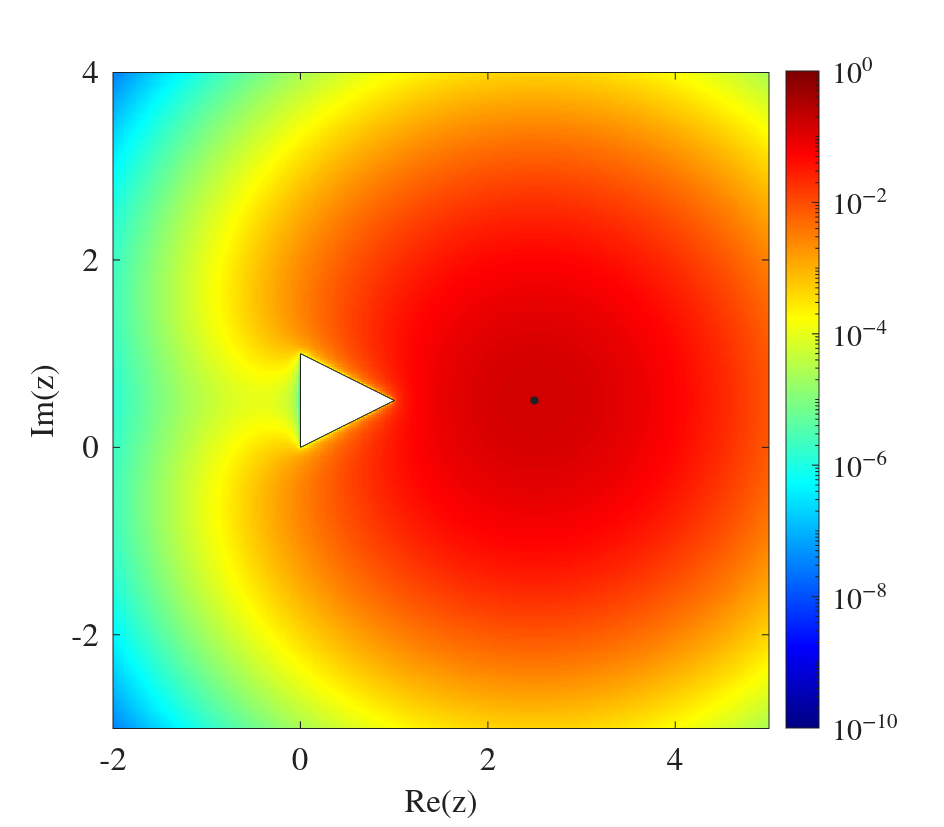}} \quad
    \subfigure[$u(z,t)$ at $t=0.8$.]{\includegraphics[width=0.31\textwidth]{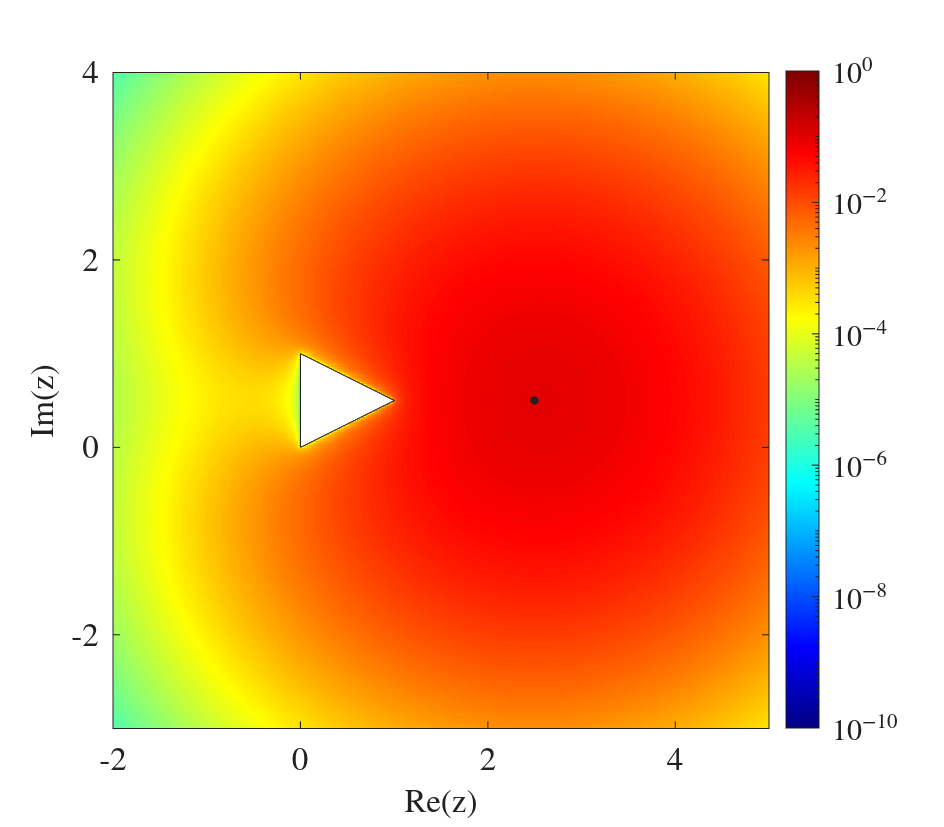}} \quad  
    \subfigure[$u(z,t)$ at $t=2$.]{\includegraphics[width=0.31\textwidth]{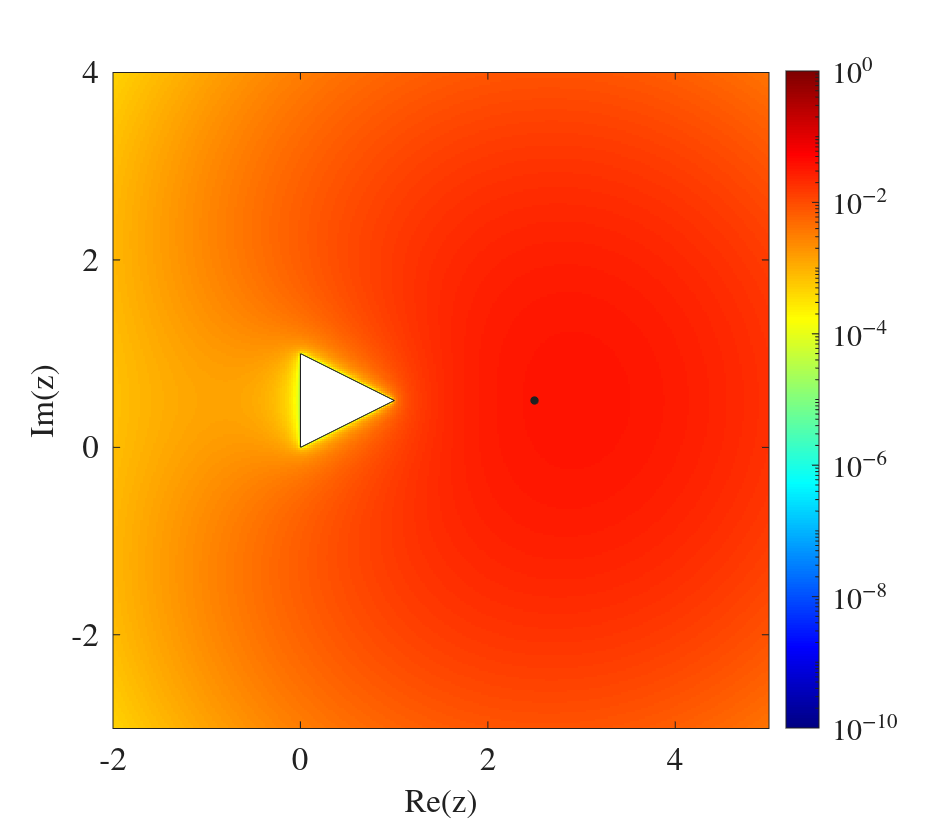}}
    \caption{LM solution $u(z,t)$ of \eqref{eq:main} for the case of a single triangular absorbing body with initial condition $u_0(z) = \delta(z-(2.5+0.5i))$, $m=90$. Solutions shown a sequence of times from $t=0.01$ to $t=2$.\label{fig:triangle} }
\end{figure}

In this example, we consider a single triangular absorber centered at $z=0.5+0.5i$ with vertices $\{1+0.5i,i,0\}$. The centering point $z_{*}$ of the Runge expansion is the centroid of the triangle. The solution $u(z,t)$ of \eqref{eq:main} is calculated for $f\equiv0$ and $u_0(z) = \delta(z-z_0)$ where $z_0=2.5+0.5i$. In Fig.~\ref{fig:triangle}, we plot solutions at $6$ time points showing shortly after initialization and the engagement with the body. In particular, we remark that the method captures the \lq\lq shadow\rq\rq\ region in the rear of the body where the solution has lower value. Probabilistically, the lower values reflects the observation that a diffusing particle is much less likely to be present at this location.

\subsection{Example: \hl{Flux calculation comparison of the LM with KMC and matched asymptotic expansion}}

In this example we compare the solution of \eqref{eq:main} from the LM with two recently developed alternative methods, a particle based Kinetic Monte Carlo (KMC) method and a matched asymptotic expansion. The asymptotic approximation, as summarized in Sec.~\ref{sec:MAA}, is derived in the limit of well separated bodies \cite{che-lin-her-qua2022}. In our experiments, we calculate and compare the cumulative density functions $c(t)$, defined in \eqref{eqn:bdflux}, with each of these methods.

\begin{figure}[h]
    \centering
    \subfigure[$u(z,t)$ at $t=0.01$.]{\includegraphics[width=0.45\textwidth]{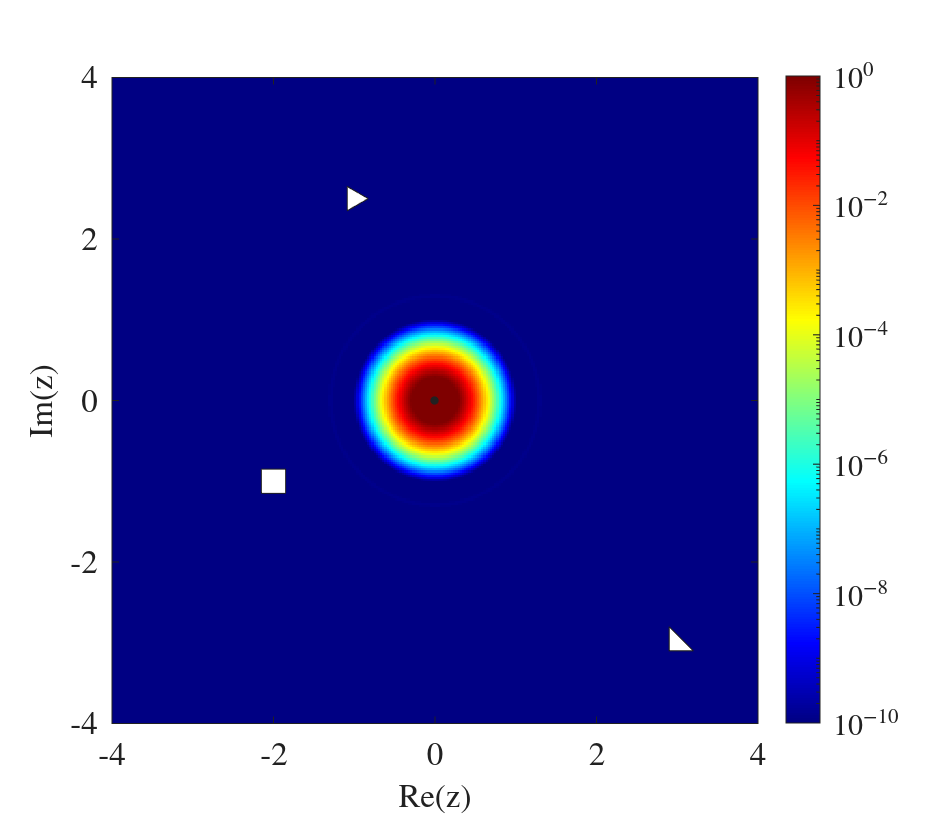}\label{fig:3_bodies_b_a}} \qquad
    \subfigure[$u(z,t)$ at $t=0.1$.]{\includegraphics[width=0.45\textwidth]{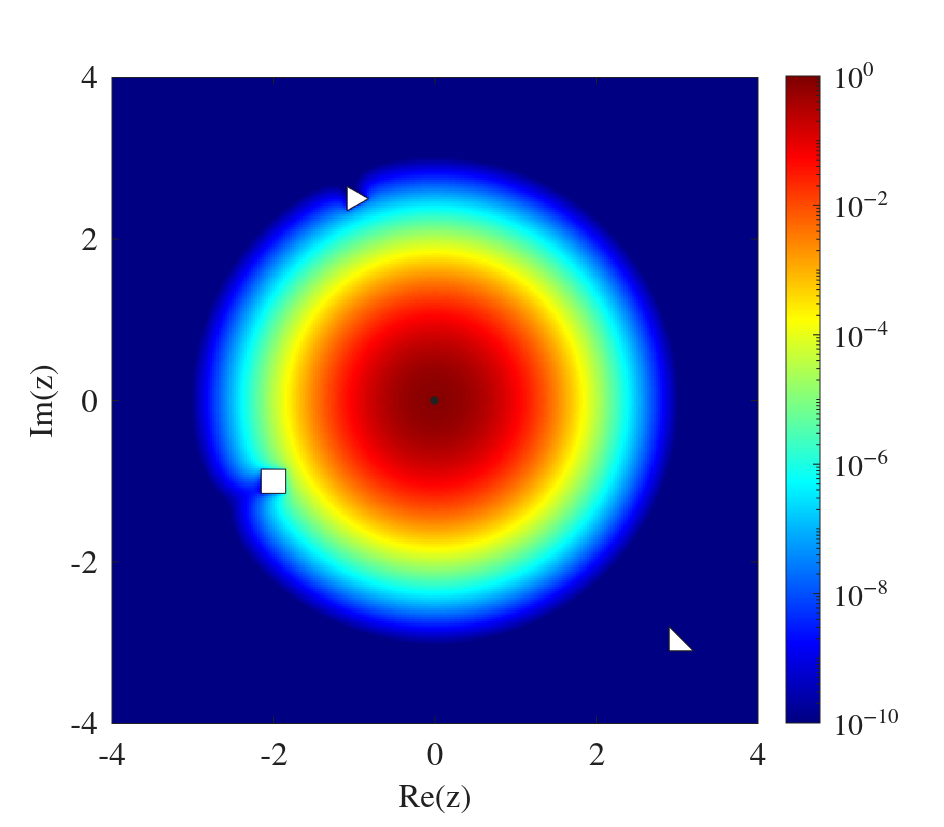}\label{fig:3_bodies_b_b}} \\
    \subfigure[$u(z,t)$ at $t=1$.]{\includegraphics[width=0.45\textwidth]{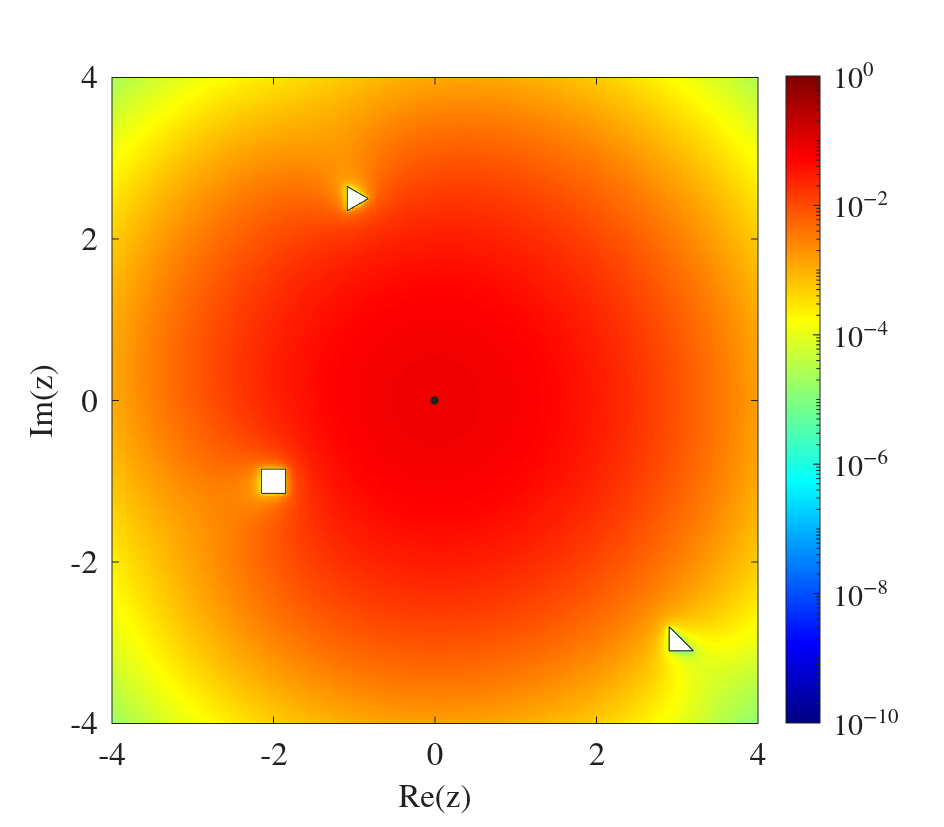}\label{fig:3_bodies_b_c}} \qquad
    \subfigure[Agreement between different methods.]{\includegraphics[width=0.45\textwidth]{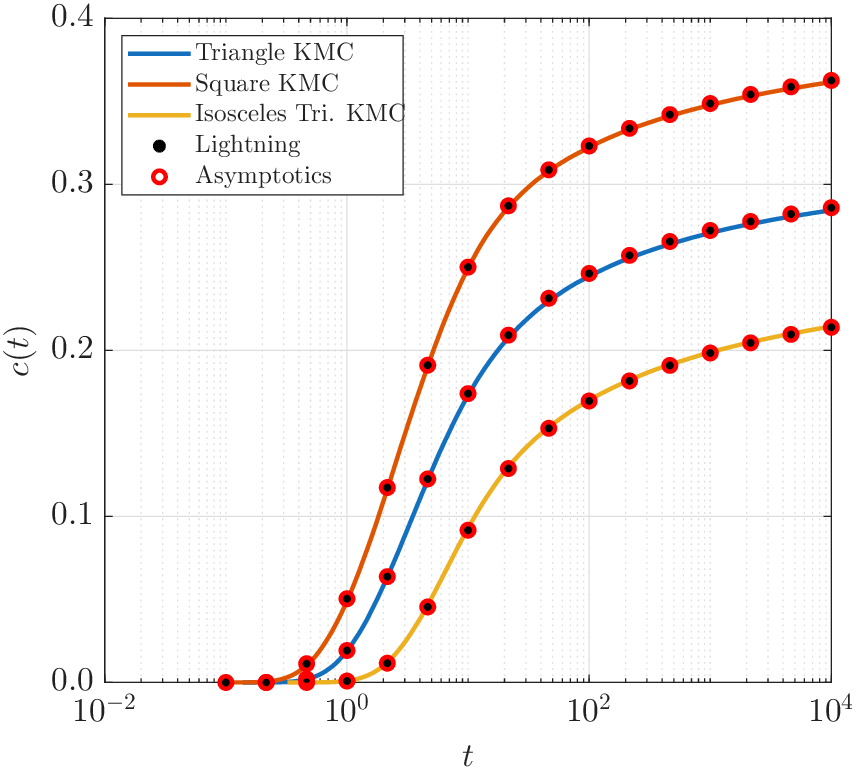}\label{fig:3_bodies_b}}
	\caption{Panels (a-c): Solutions $u(z,t)$ of \eqref{eq:main} with three small absorbers with $u_0(z)=\delta(z)$ and \hl{$m=60$} at various $t$ from $0.01$ to $1$. Panel (d): Agreement of the capture rate \eqref{eqn:bdflux} at individual absorbers as computed from the LM, the asymptotic solution \eqref{eqn:asy_solution} and the KMC method \cite{che-lin-her-qua2022}.\label{fig:3_bodies}}
\end{figure}

Our first configuration consists of a point source at $z_0=0$ diffusing to $N_B=3$ bodies; an equilateral triangle centered at $z_1-1+2.5i$, a square centered at $z_2 = -2-i$, and an isosceles triangle with centroid at $z_3 = 3-3i$. The side lengths of the square and triangle, and the small-side length of the isosceles triangle are uniformly $h=0.3$. The flux vector is calculated from the asymptotic formula \eqref{eqn:asy_solution} and then inverted with $M=9$ points. The KMC simulation is run with $1\times 10^6$ particles. In Fig.~\ref{fig:3_bodies}, we plot the solution $u(z,t)$ at $t=1$ and again highlight the shadow effect where the occupation probability behind the target is markedly lower. The agreement between the cumulative fluxes $c(t)$ from the LM, KMC and asymptotic methods is very good (cf.~Fig.~\ref{fig:3_bodies_b}). 

\begin{figure}[h]
    \centering
    \subfigure[$u(z,t)$ at $t=0.05$.]{\includegraphics[width=0.45\textwidth]{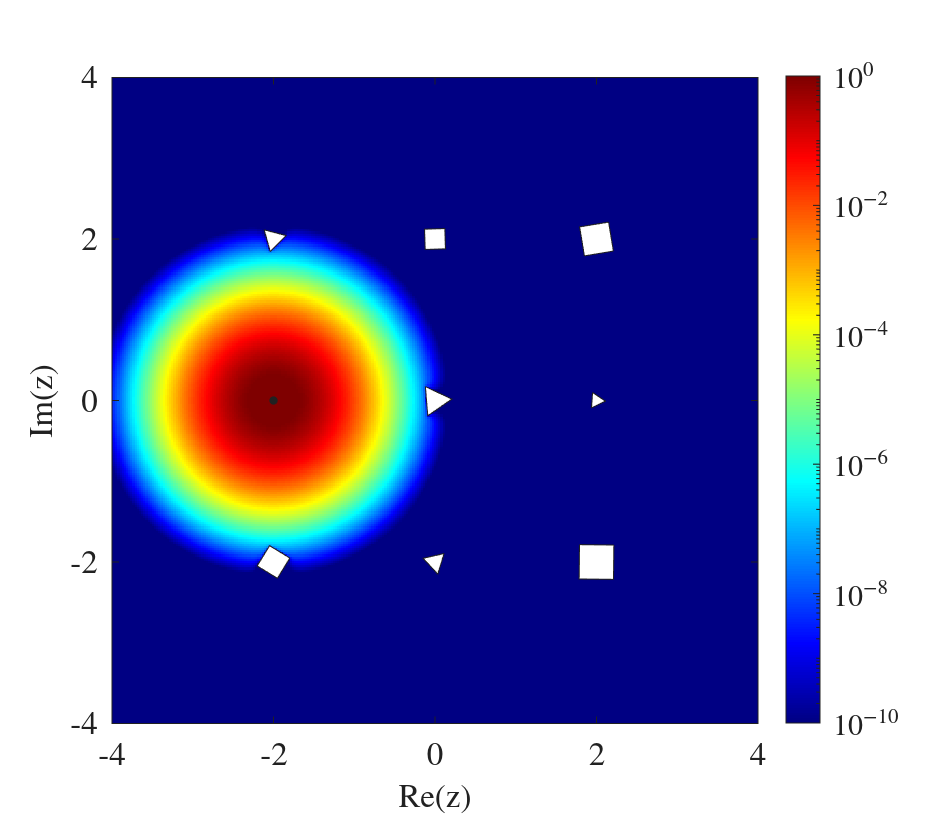}}\qquad
    \subfigure[$u(z,t)$ at $t=0.2$.]{\includegraphics[width=0.45\textwidth]{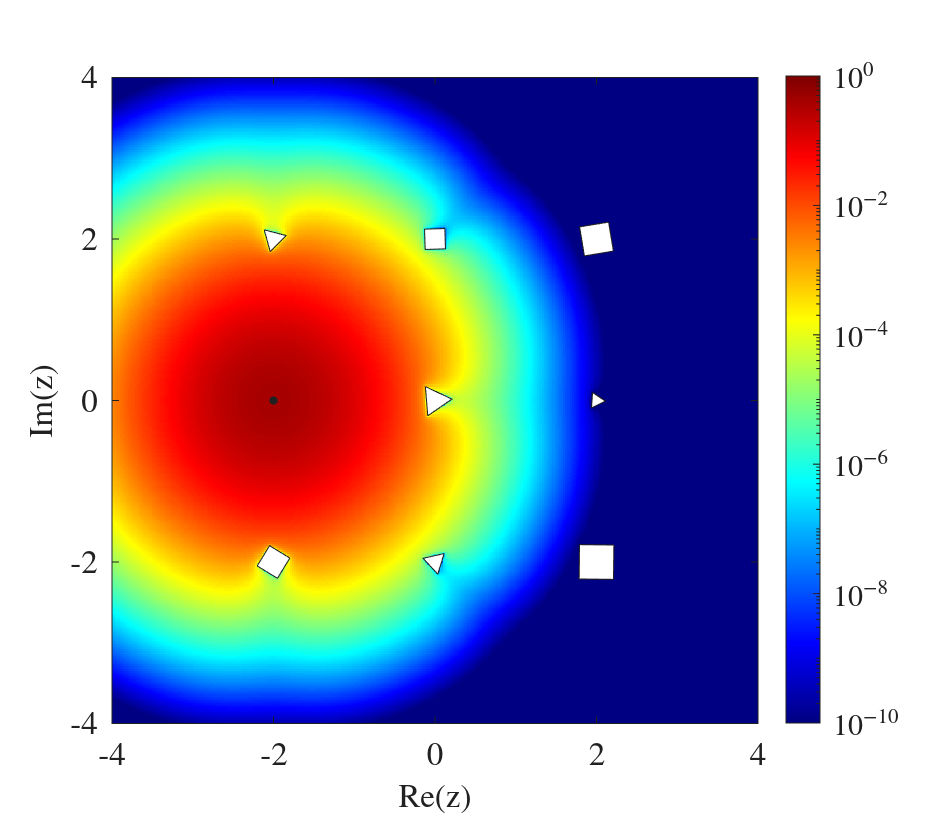}}\\
    \subfigure[$u(z,t)$ at $t=1$.]{\includegraphics[width=0.45\textwidth]{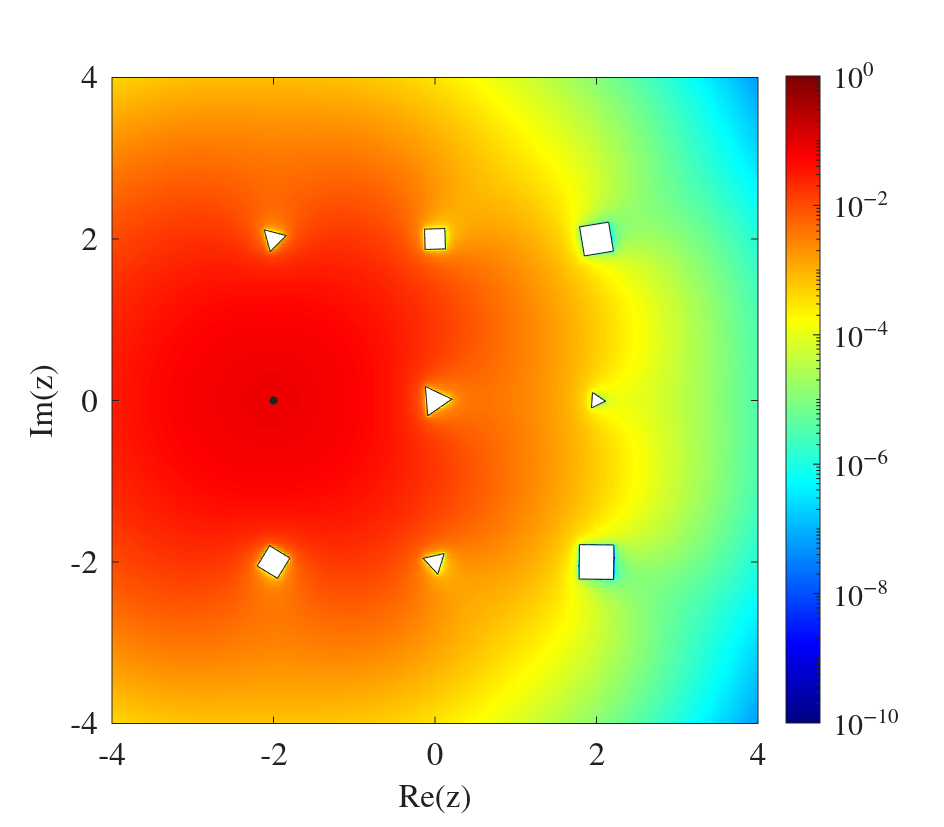}}\qquad
    \subfigure[Agreement between LM, KMC and asymptotics.]{\includegraphics[width=0.45\textwidth]{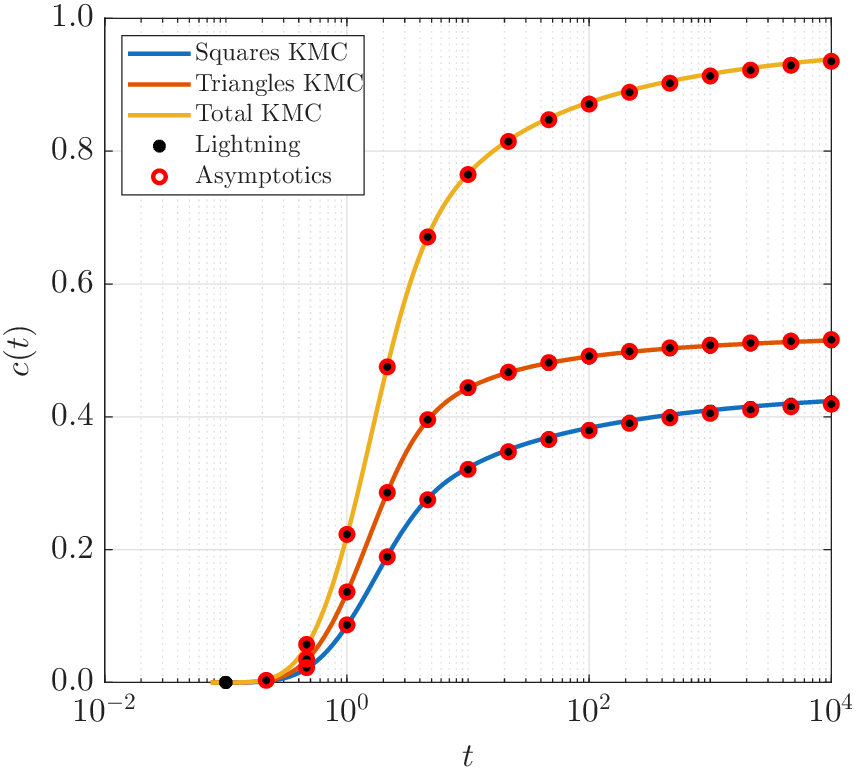} \label{fig:8_bodies_d}}
	\caption{Panels (a-c): LM solution $u(z,t)$ of \eqref{eq:main} with \hl{$m=60$}. Panel (d): Cumulative fluxes $c(t)$ through all squares, all triangles, and all bodies computed by LM, KMC and matched asymptotic expansions.\label{fig:8_bodies}}
\end{figure}

For our next configuration, we consider $N_B=8$ bodies of various $h$ and rotations in an integer grid with spacing $2$ and a point source at $z_0 = -2$. Rotations are sampled uniformly on $[0,2 \pi)$ and side lengths are normally distributed with mean $0.3$ and standard deviation $0.2$ with positivity enforced. The domain and three solutions $u(z,t)$ of \eqref{eq:main} are shown in Figs.~\ref{fig:8_bodies}(a-c). Solutions are also calculated from the matched asymptotic expansion and the KMC method. We choose half of the bodies to be equilateral triangles, and half to be squares, and calculate their combined fluxes, once again attaining very close agreement (cf.~Fig.~\ref{fig:8_bodies_d}). Given the fact that these three solutions are generated by completely different methods, their remarkably good agreement serves as a highly non-trivial validation of the LM approach \hl{for well-separated small bodies}.

\begin{figure}[h]
    \centering
    \subfigure[$\phi_1(z_0)=0.2552$]{\includegraphics[width=0.45\textwidth]{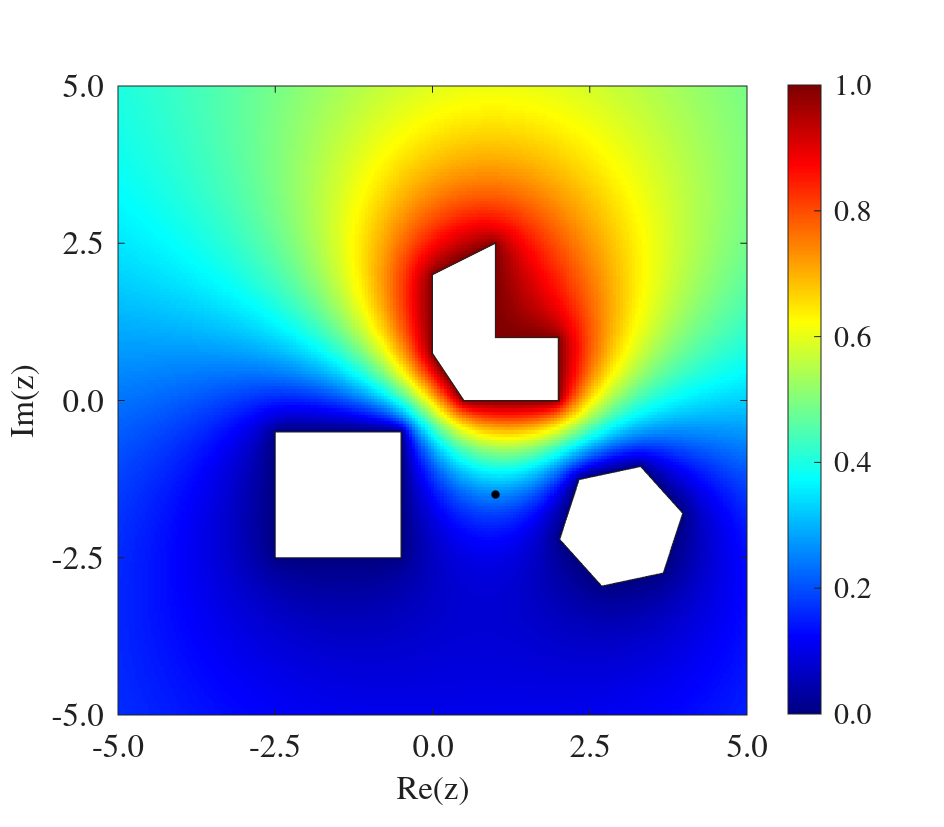}\label{fig:split_a}} \qquad
    \subfigure[$\phi_2(z_0)=0.3439$]{\includegraphics[width=0.45\textwidth]{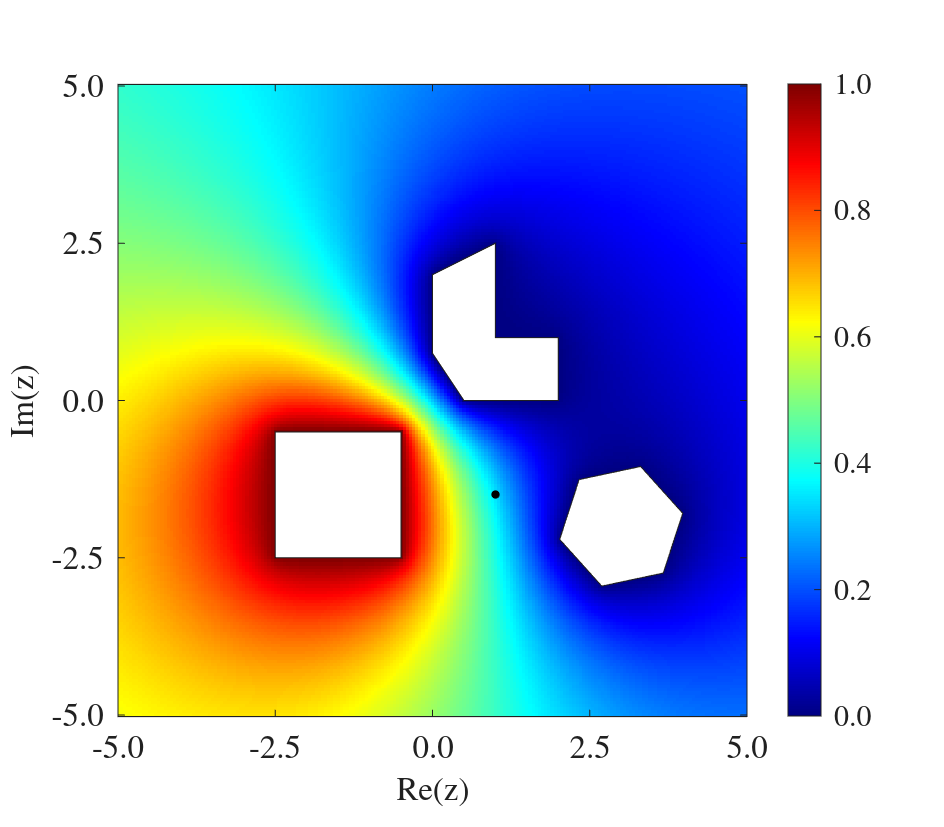}\label{fig:split_b}} \\
    \subfigure[$\phi_3(z_0)=0.4009$]{\includegraphics[width=0.45\textwidth]{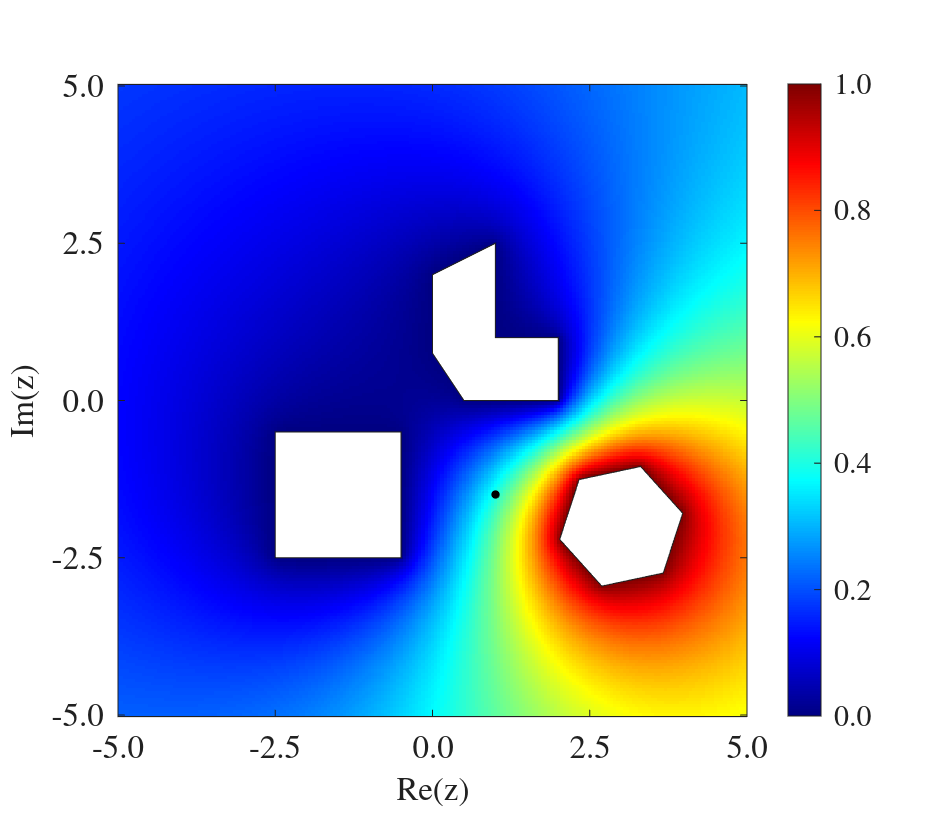}\label{fig:split_c}} \qquad
    \subfigure[Agreement between LM and KMC]{\includegraphics[width=0.45\textwidth]{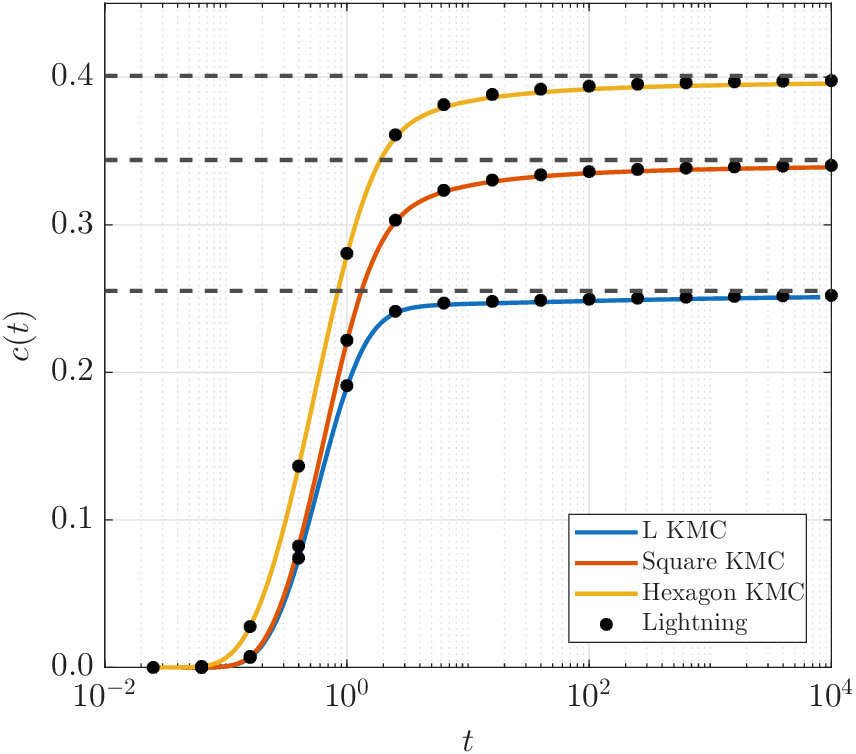}\label{fig:splits}}
	\caption{\al{Panels (a-c): Solutions $\mathrm{Re}(\phi_j(z))$ of the splitting problems \eqref{eq_splitting}  for the L shape, square, and hexagon respectively solved via the Lightning Laplace method with $m=60$. The black dot denotes the point source $z_0$ for panel (d): The cumulative fluxes calculated from KMC and the lightning method solution. Dotted horizontal lines are drawn at the splitting probabilities evaluated at $z_0$. In particular, we note the agreement between $\lim _{t \to \infty} c(t)$ and $\phi_j(z_0)$ providing a highly non-trivial validation of our method.\label{Fig:Splitting}}}
\end{figure}

\al{As an additional validation, we calculate the so-called splitting probabilities $\phi_j(z)$, for $j= 1,\ldots,N_B$. These quantities determine the likelihood \cite{Venu2015,Lindsay2023a} that a diffusing particle initially at $z\in \mathbb{C}\setminus\Omega$ will eventually be absorbed at the $j^{\text{th}}$ body. They satisfy the Laplace problems
\bsub\label{eq_splitting}
\begin{gather}
\label{eq_splitting_a} \Delta \phi_j = 0, \quad z\in \mathbb{C}\setminus\Omega; \qquad \phi_j (z) \text{ bounded as } |z|\to\infty;   \\[5pt]
\label{eq_splitting_b} \phi_j= \delta_{j,k}, \quad z\in\partial\Omega_k; \qquad k = 1,\ldots, N_B.
\end{gather}
\esub}
\hl{In Fig.~\ref{Fig:Splitting} we plot the solutions of \eqref{eq_splitting} for a new case of $N_B=3$ bodies where the bodies are neither small or well separated and thus not suitable for calculating $c(t)$ by matched asymptotic expansion. These solutions are acquired via the lightning method for the Laplace problem, in which pole prescriptions remain the same, but $\psi_k$ are instead rational terms. One must also caution that the exterior domain Laplace problem necessitates a logarithmic potential in our Runge expansion (see \cite{cos-tre2023}, which solves it with the related AAA method). Also plotted in the figure is the cumulative fluxes captured by these bodies with a point source at $z_0=1-1.5i$. We observe that the limiting behavior of the cumulative flux into each body asymptotically approaches the splitting probabilities evaluated as $t\to\infty$. This agreement serves as a highly non-trivial verification of the entire solution approach including the evaluation of integrated fluxes.}

\subsection{Example: Nonzero boundary conditions.}\label{sec:non_zero}

\begin{figure}[htbp]
    \centering
    \subfigure[$|\mathrm{Re}(\hat{u}(z,s))|$ at $s=-30-21i$.]{\includegraphics[width=0.325\textwidth]{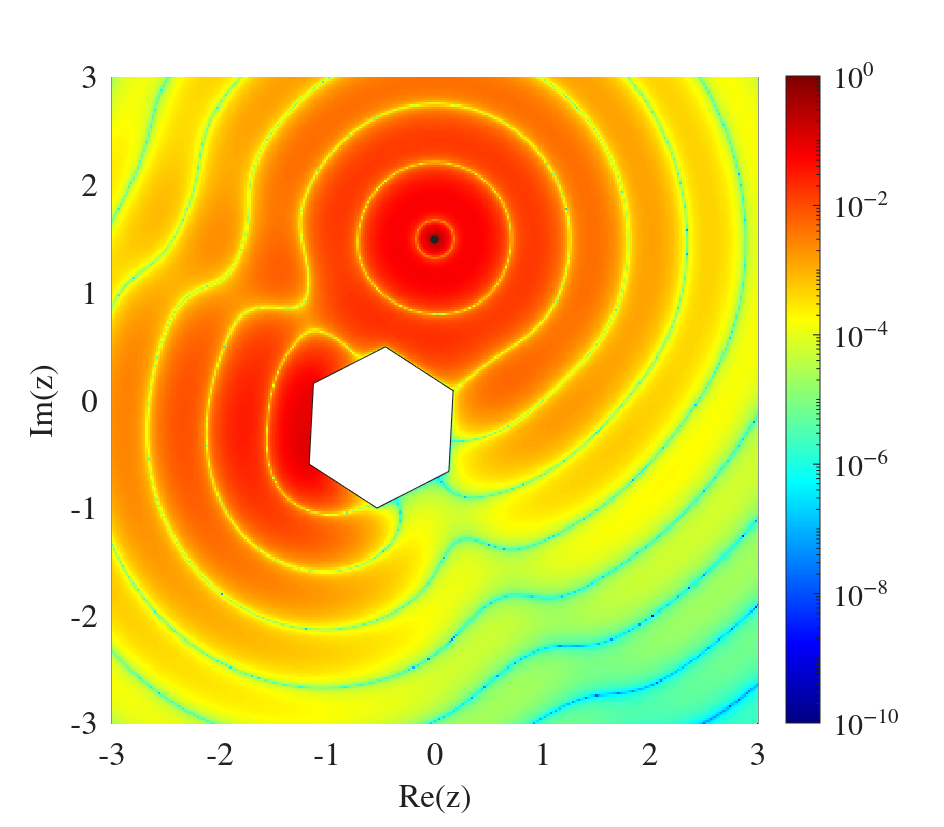}\label{fig:ex_nonzeroBC_a}} 
    \subfigure[$u(z,t)$ at $t=0.03$.]{\includegraphics[width=0.325\textwidth]{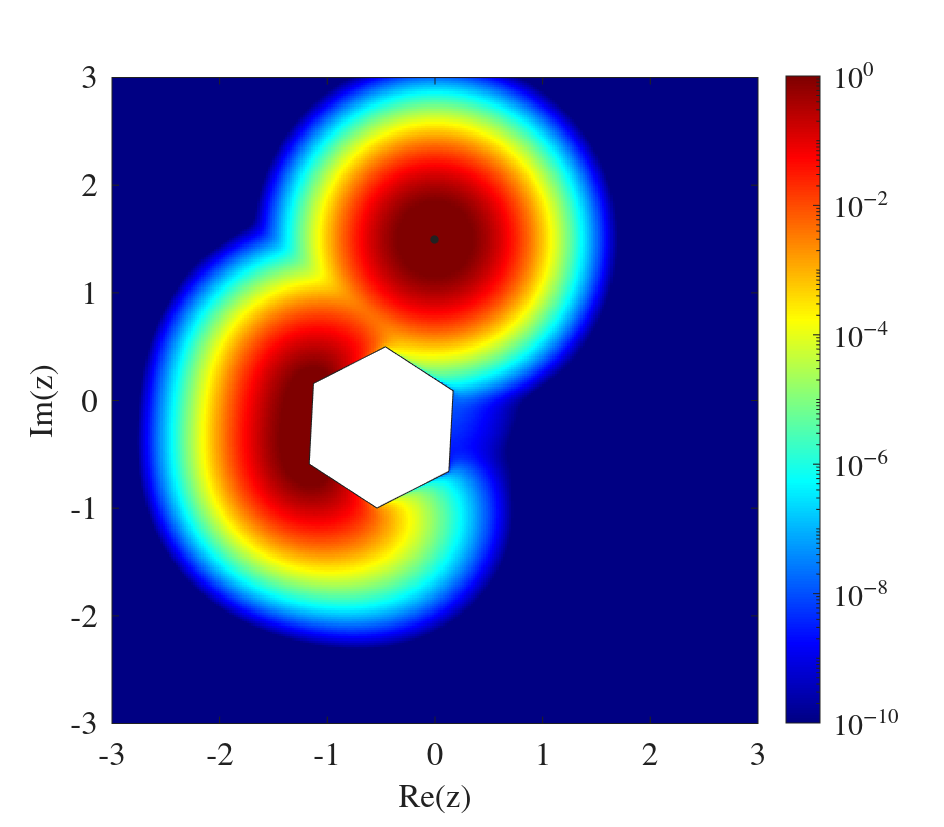}} 
    \subfigure[$u(z,t)$ at $t=0.1$.]{\includegraphics[width=0.325\textwidth]{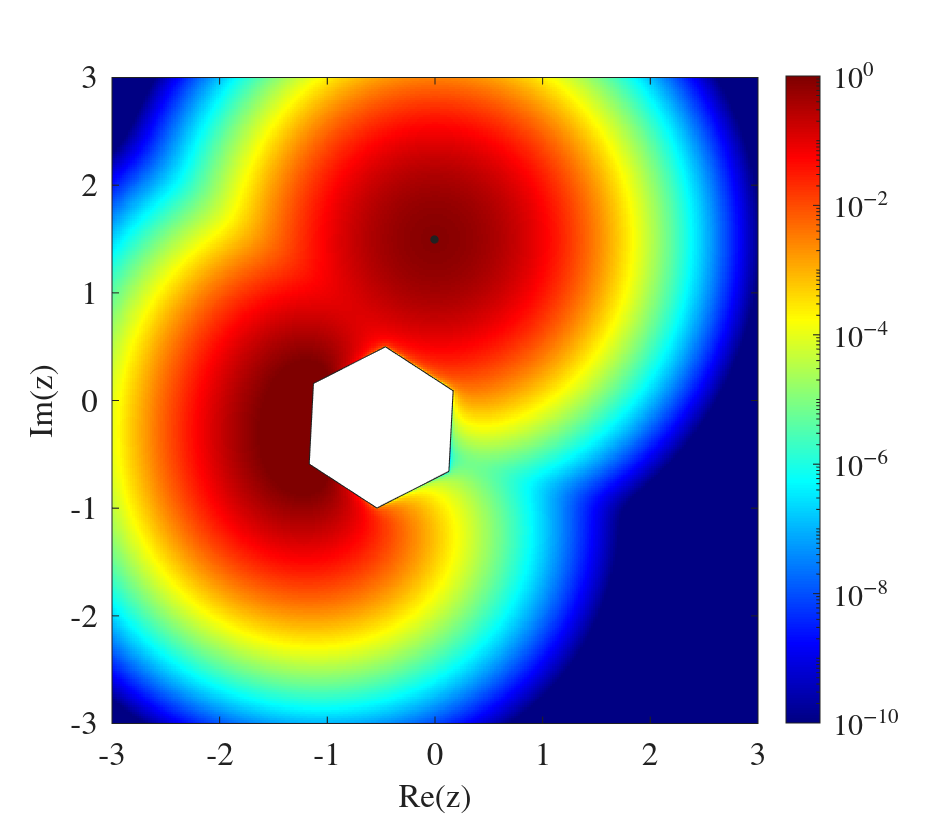}}
    \caption{The LM method for the solution of \eqref{eq:main} with a single hexagonal body, the applied boundary condition $u = f(z)=\mathrm{Re}(z)^{10}$ on $\partial\Omega$ and the initial condition $u_0 = \delta(z-z_0)$ for $z_0 =1.5i $ and $m=90$. Panel (a): The solution of \eqref{eq:Helmholtz} at value $s=-30-21i$. Panel (b): The solution of the heat equation \eqref{eq:main} at $t=0.03$. Panel (c) The solution $u(z,t)$ at $t=0.1$. \label{fig:ex_nonzeroBC} }
\end{figure}

Our previous examples focused on the homogeneous boundary condition $u(z,t) = f(z) \equiv 0$ on $\partial\Omega$ and the initial condition $u(z,0) = \delta(z-z_0)$. In this section, we demonstrate the application of the LM method to nonzero boundary conditions. We impose the condition $f(z)=\mathrm{Re}(z)^{10}$ In Fig.~\ref{fig:ex_nonzeroBC} we plot the solution $u(z,t)$ of \eqref{eq:main} at several time points for the initial condition $u_0(z) = \delta(z-z_0)$ and $z_0 = 1.5$.

\subsection{Example: L-shaped domain and multiple Runge poles.}
As remarked in other works \cite{Gopal19B}, elongated domains pose a difficulty to lightning solvers due to the Runge expansion decaying as it approaches $\partial \Omega$ and becoming an ill-conditioned basis. As $K_n(z)$ decays exponentially as $z \to \infty$, rather than the polynomial rate of rational functions, \hl{we observe increased ill-conditioning}.

\begin{figure}[htbp] 
    \centering    
    \subfigure[$u(z,t)$ at $t=0.03$ and $z_0=1.5+1.5i$.]{\includegraphics[width=0.325\textwidth]{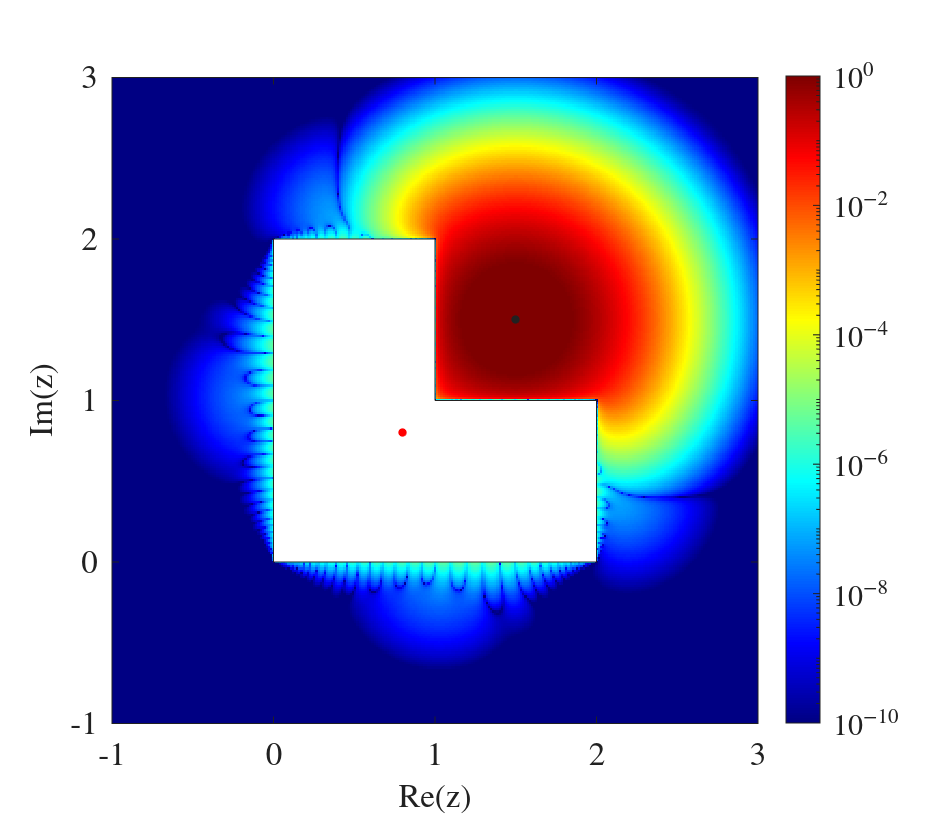}}
    \subfigure[$u(z,t)$ at $t=0.03$ and $z_0 = 2+2i$.] {\includegraphics[width=0.325\textwidth]{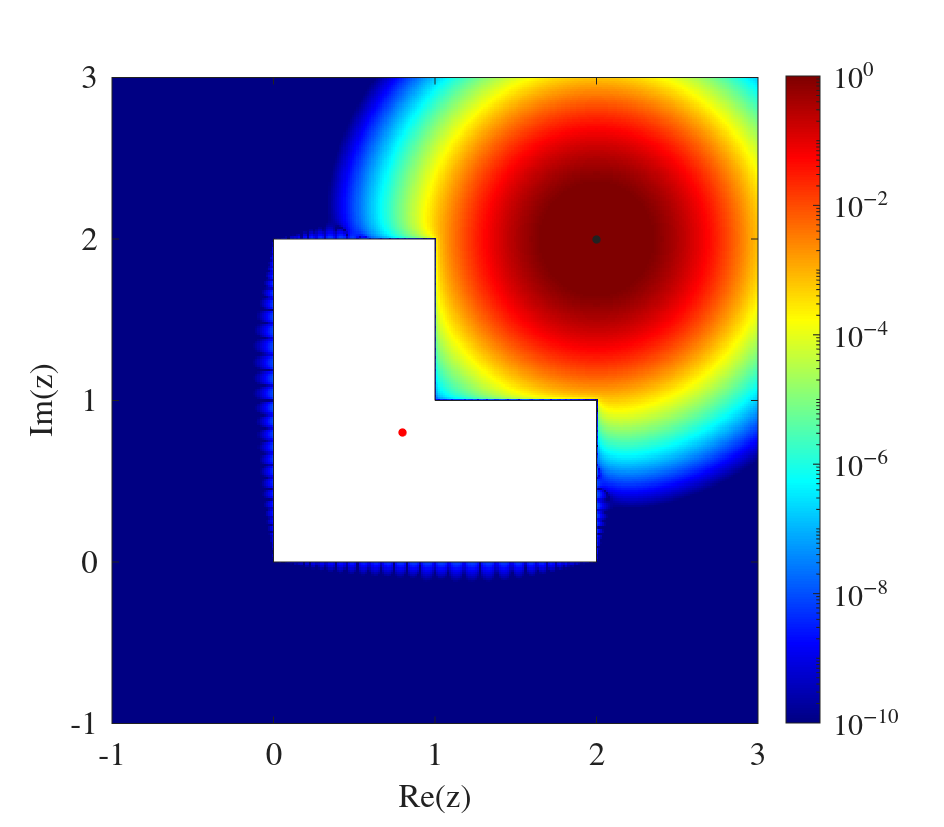}}
    \subfigure[$u(z,t)$ at $t=5$ and $z_0 = 1.5+1.5i$.]{\includegraphics[width=0.325\textwidth]{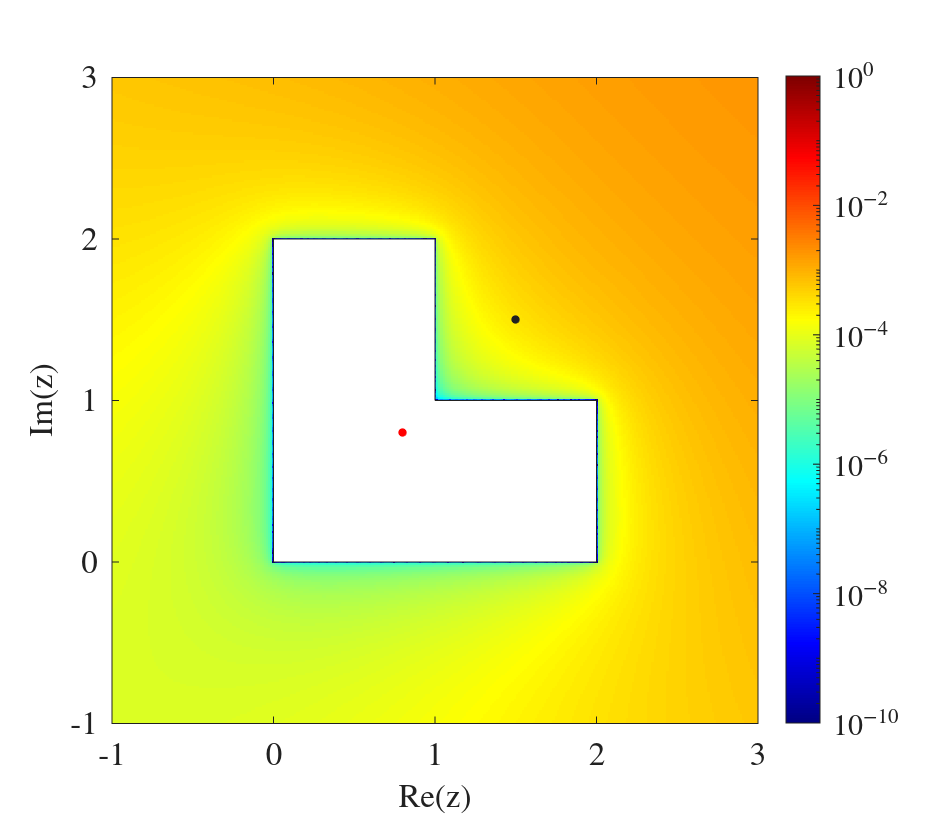}}    
	\caption{Solutions of \eqref{eq:main} for the L-shaped domain with $m=90$, $u_0 = \delta(z-z_0)$ for various $z_0$ (solid black dot) and different times $t$. Calculations performed with a single Runge expansion centered at $z_*=0.8+0.8i$ (solid red dot). The resulting relative errors are \hl{(a) $\mathcal{E}^{\infty}[u] = 1.72\times10^{-5}$, (b) $\einf[u] = 2.6\times 10^{-5}$, (c) $\einf[u] =6.73 \times 10^{-6}$.} \label{fig:Lsinglerunge}} 
\end{figure}

This phenomenon is demonstrated in the convergence rate for the case of the L-shaped domain with a single Runge pole at $z_*=0.8+0.8i$ (cf.~Fig.~\ref{fig:Lsinglerunge}). We solve this problem at various times and for configurations of source locations, both close and far from the re-entrant corner. We observe in each case, the root exponential convergence is maintained, albeit at a much slower rate and we do not reach our target error $10^{-10}$ \hl{before the method stagnates}, as shown in Fig.~\ref{fig:Lshape_cvg}.

\begin{figure}[htbp]
    \centering    
    \subfigure[$u(z,t)$ at $t=0.03$ and $z_0=1.5+1.5i$.]{\includegraphics[width=0.325\textwidth]{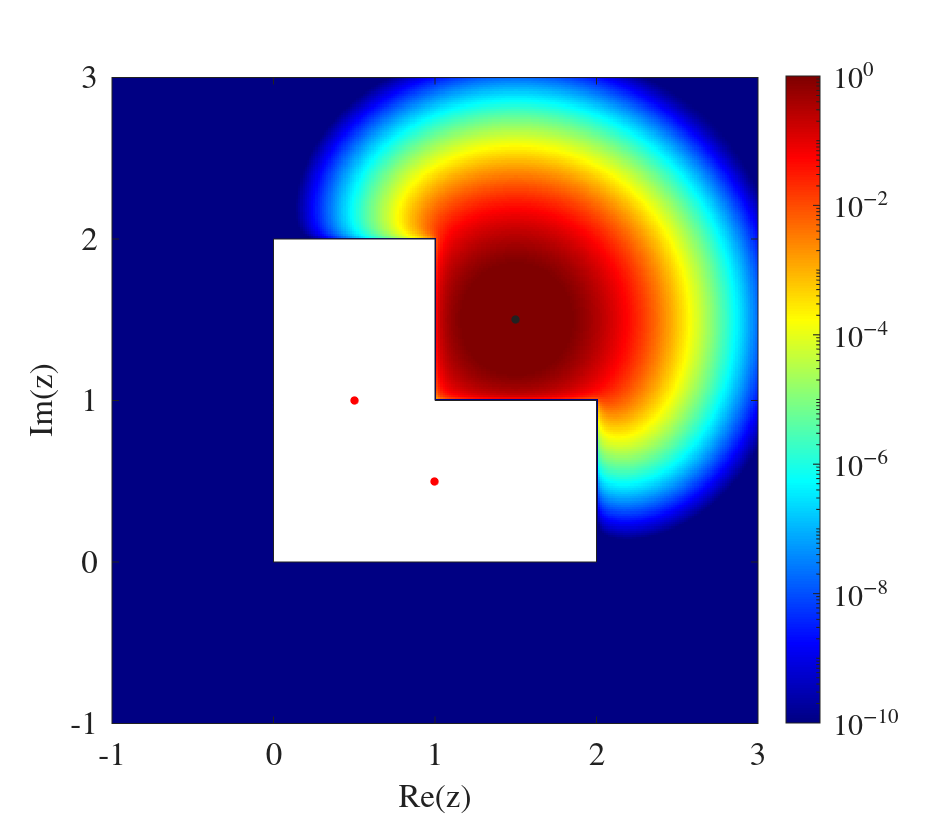}}
    \subfigure[$u(z,t)$ at $t=0.03$ and $z_0=2+2i$.]{\includegraphics[width=0.325\textwidth]{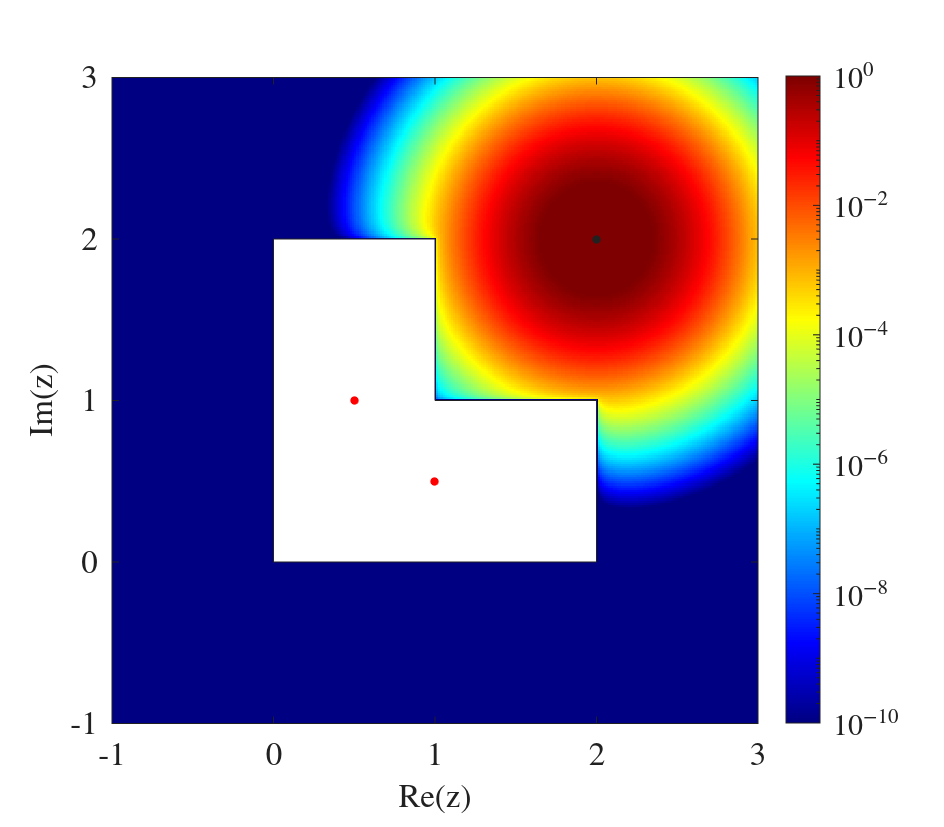}}
    \subfigure[$u(z,t)$ at $t=5$ and $z_0=1.5+1.5i$.]{\includegraphics[width=0.325\textwidth]{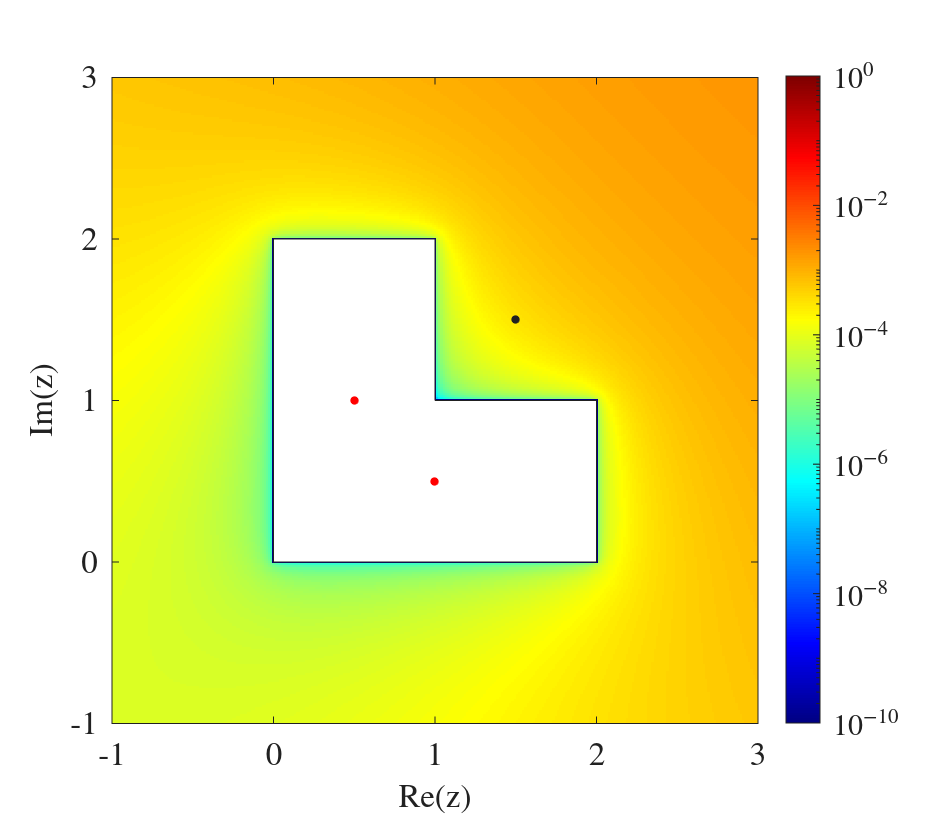}}
    
	\caption{The L-shaped domain with two Runge expansions about points $z_{*1}=1+0.5i$ and $z_{*2}=0.5+i$ (solid red dots) with $m=90$. The resulting relative errors are \hl{ (a) $\mathcal{E}^{\infty}[u] = 4.2\times10^{-10}$, (b) $\einf[u] = 2.6\times 10^{-9}$ and (c) $\einf[u] =1.36 \times 10^{-10}$. }\label{fig:Ldoublerunge}} 
\end{figure}

\begin{figure}[htbp]
    \centering    
    \includegraphics[width=0.6\textwidth]{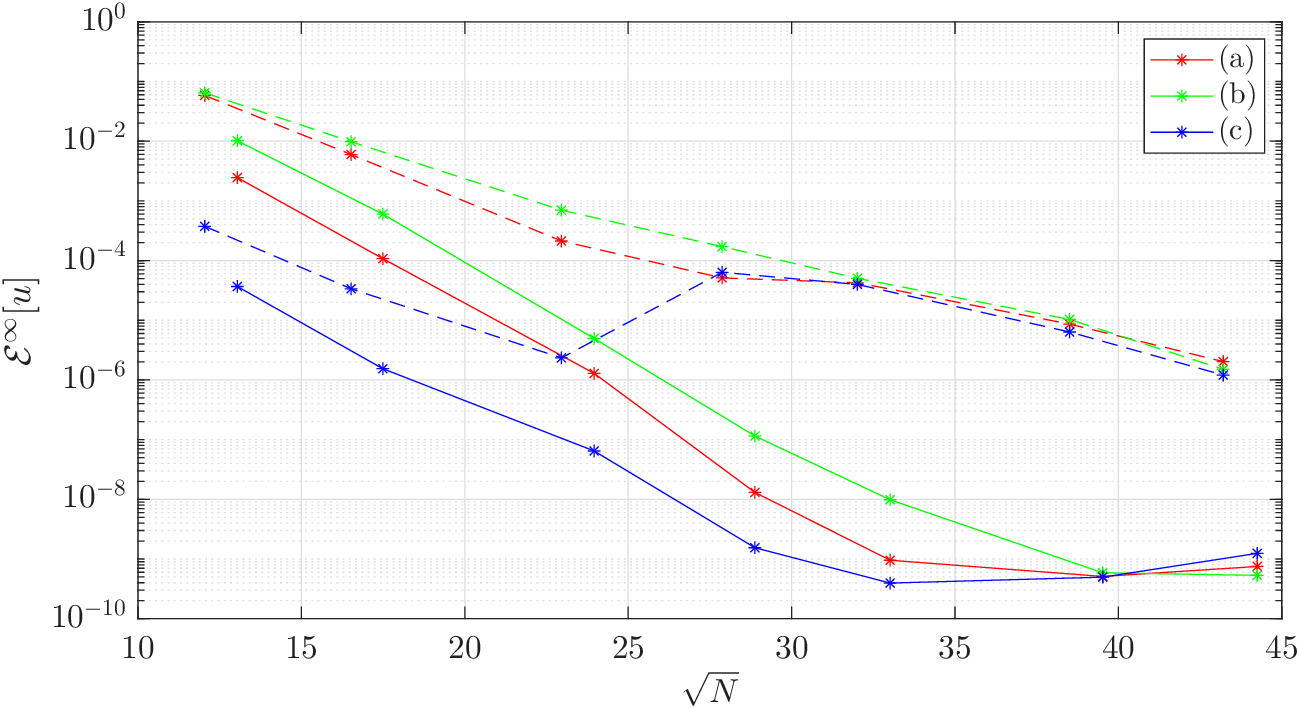}
	\caption{\hl{Convergence studies on the L shaped domain for the single Runge expansion solution (dotted lines) and the two Runge expansion solution (solid lines), on the configurations in Figs.~\ref{fig:Lsinglerunge} and \ref{fig:Ldoublerunge}. Attained convergence rate is root exponential, albeit at slower rates than the solution with a single Runge expansion.\label{fig:Lshape_cvg}}} 
\end{figure}

\hl{In the case of the Laplace problem, the Runge part of the problem matrix is in the form of a Vandermonde matrix, whose conditioning can be improved by Arnoldi orthogonalization \cite{brubeck2021vandermonde}. For the solution of \eqref{eq:Helmholtz}, it may be possible to modify this orthogonalization procedure to be compatible with our basis of Bessel functions. In lieu of this approach, we obtain a simple resolution of this issue by placing multiple Runge expansions into our series representation of the solution. In Fig.~\ref{fig:Ldoublerunge}, we show results for two Runge expansions centered at the points $z_{*1} = 1+0.5i$ and $z_{*2} = 0.5 + i$ and observe an improved rate of convergence and reduction in error.}. As seen in the solution profiles $u(z,t)$ in Figs.~\ref{fig:Ldoublerunge}(a-c), the solution is now well resolved near the boundary. Indeed, in Fig.~\ref{fig:Lshape_cvg}, we see a dramatic improvement in the convergence rate of the method.

\subsection{Example: Solution for general initial conditions.}
\begin{figure}[htbp]
    \centering
    \subfigure[Initial profile $u_0(z)$ and domain.]{\includegraphics[width=0.46\textwidth]{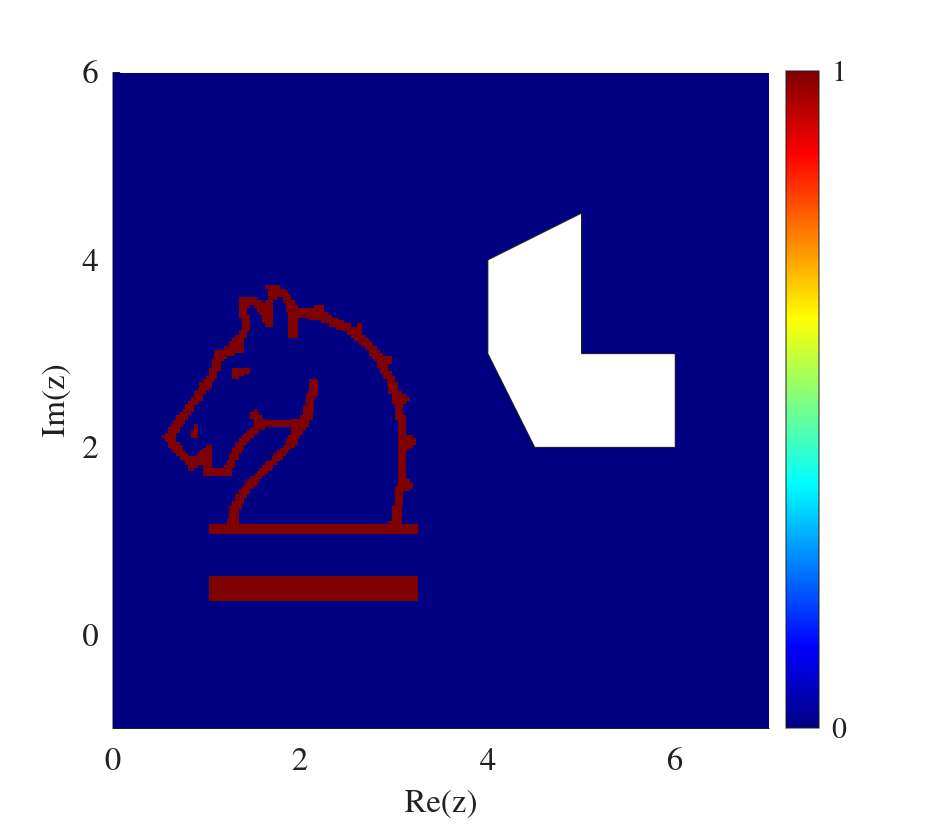}\label{fig:u0_za}} \qquad
    \subfigure[$|\mathrm{Re}(\hat{u}(z;s))|$ for $s=-68+249i$.]{\includegraphics[width=0.49\textwidth]{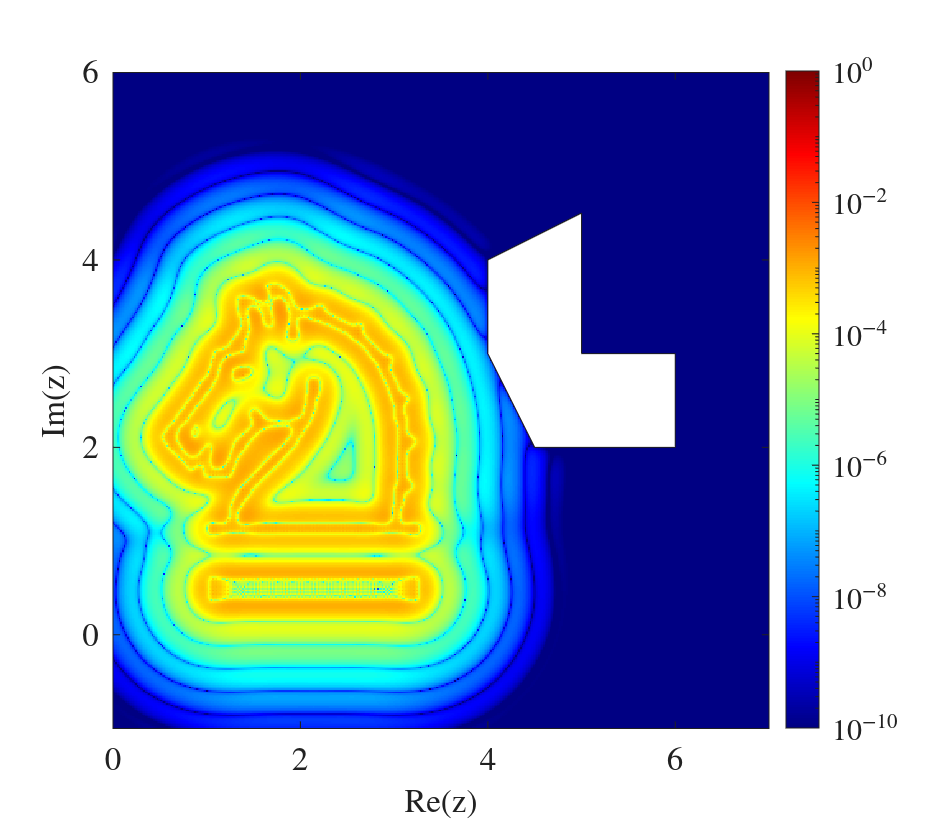}} \\
    \subfigure[$u(z,t)$ at $t = 0.03$.]{\includegraphics[width=0.325\textwidth]{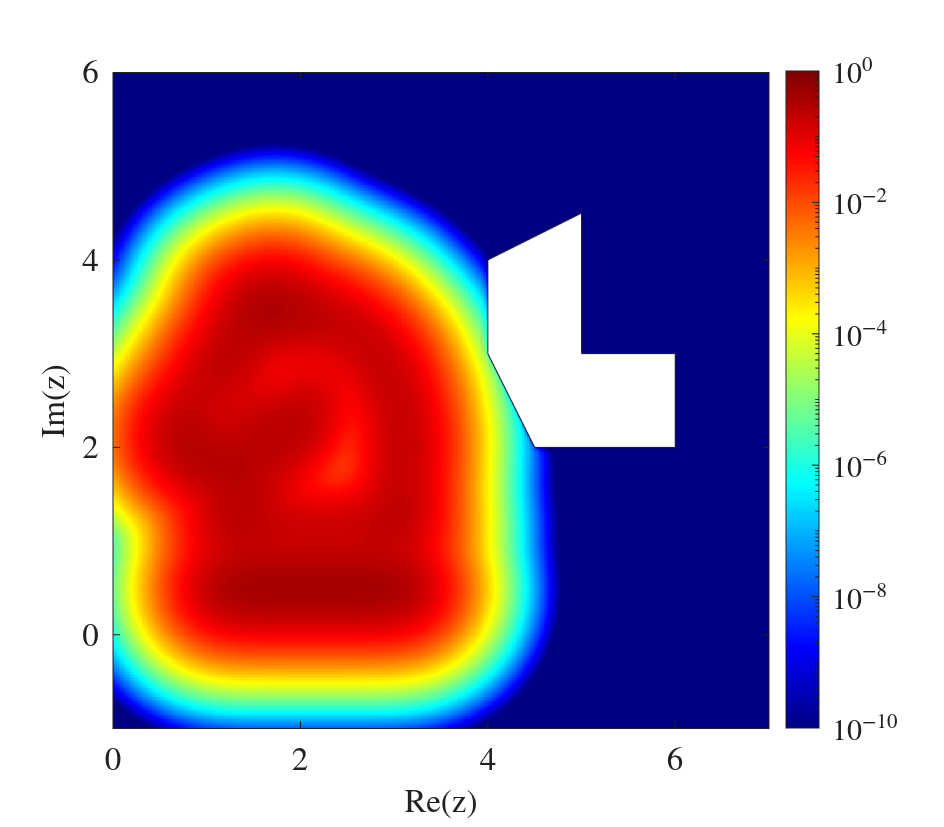}}
    \subfigure[$u(z,t)$ at $t = 0.1$.]{\includegraphics[width=0.325\textwidth]{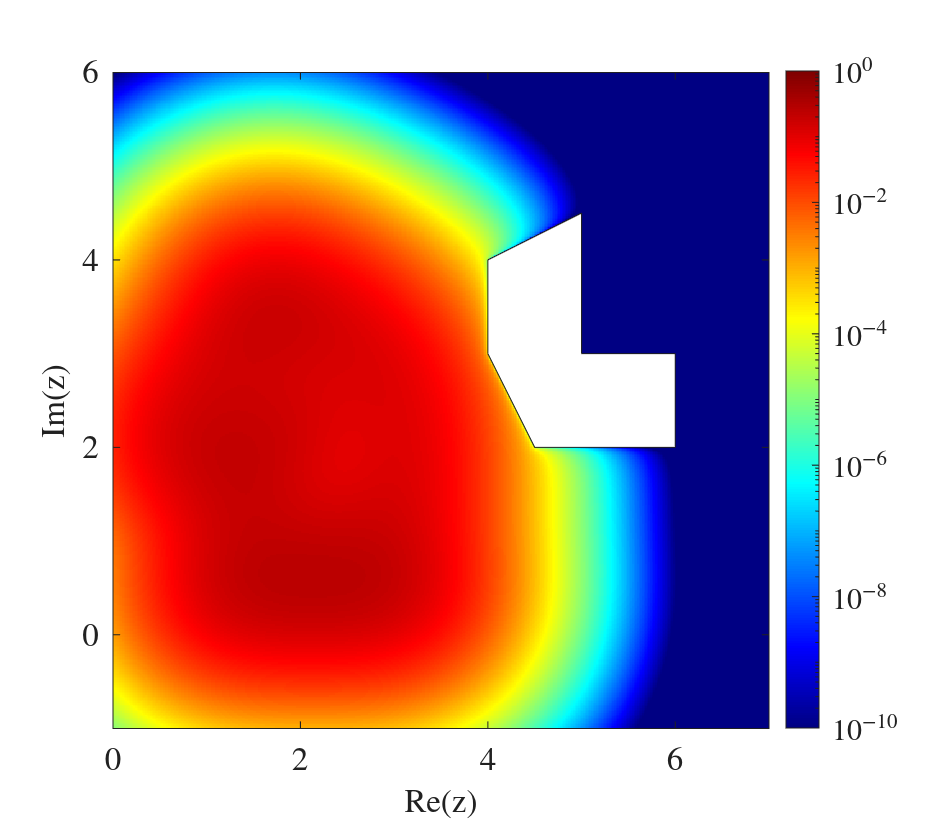}}
    \subfigure[$u(z,t)$ at $t = 1$.]{\includegraphics[width=0.325\textwidth]{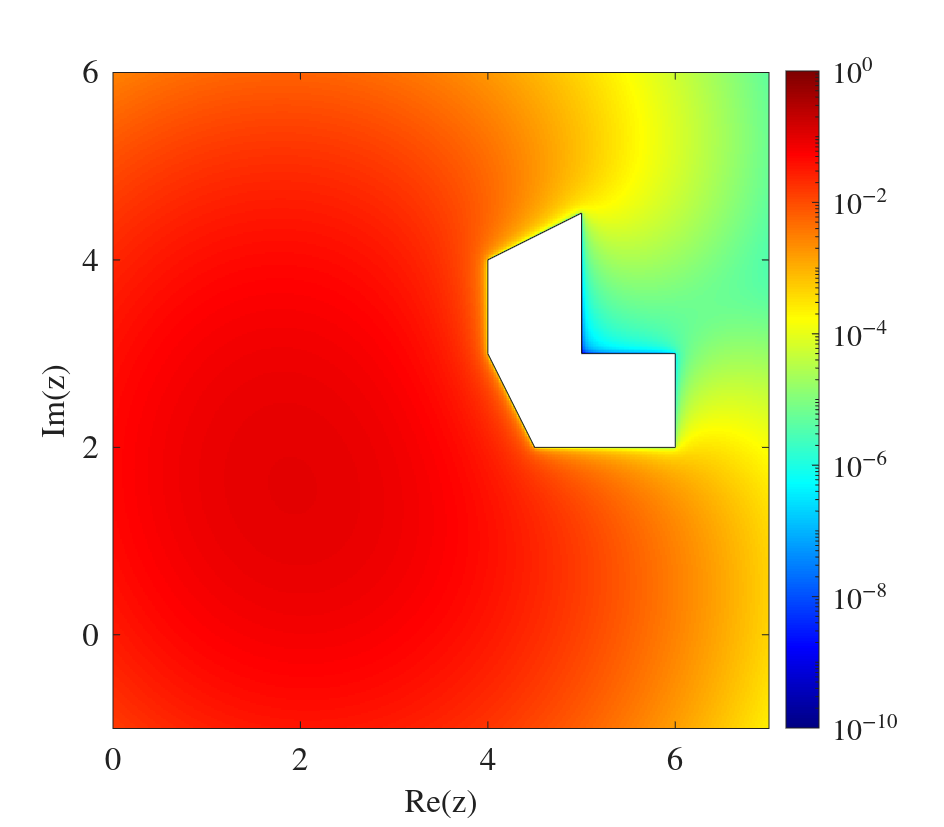}}
    \caption{The LM applied to equation \eqref{eq:main} for $f=0$, a single polygonal body and a general initial condition, all solved with $m=90$. Panel (a): The domain and the initial profile $u_0(z)$ \al{given by the characteristic function of the region}. \al{The two Runge poles are shown (solid red dots)}. Panel (b): Solution of \eqref{eq:Helmholtz} for $s=-68+249i$. Panels (c-e): solutions $u(z,t)$ of \eqref{eq:main} for various $t$. \label{fig:u0_z}}
\end{figure}

In this example, we demonstrate the applicability of the method for the more general initial condition $u_0(z)$ \al{given by the characteristic function of the region} shown in (cf.~Fig.~\ref{fig:u0_za}). This function is compactly supported on $\mathbb{R}^2$ and we evaluate the integral term $\hat{u}_p$ in \eqref{eq:particularsol} via quadrature. The shown solution profiles $u(z,t)$ display the smoothing effects of the heat equation on the discontinuities in the initial profile. For our geometry, we adopt our configuration in Fig. \ref{fig:sigma} with Runge expansions at $z_{*1}=5+2.5i$ and $z_{*2}=4.5+3i$.

\section{Discussion.}\label{dec:discussion}

In this work we have developed a lightning solver for the heat equation in unbounded planar regions with disjoint absorbing bodies. Our approach provides a complement to recent adaptations of the lightning method to Laplace, Helmholtz and Stokes problems and further highlights its usefulness in the numerical computation of solutions to an important class of challenging problems where geometric corners induce solution singularities. Our method has been validated on several complex examples and through comparison with solutions via complementary techniques such as \al{Monte Carlo methods, matched asymptotic expansions and boundary integral equations.}

There are various avenues for future work that emanate from this study. Our applications have focused solely on the case of Dirichlet absorbers, however, it is natural to consider additional boundary conditions, such as Neumann and Robin which appear in numerous biological signaling problems \cite{ye2024,Lawley2024}. In rudimentary investigations on the applicability of this method to the case of all Neumann bodies, we observe that jump discontinuities of the normal derivative at corners can induce oscillations in the solution. A similar situation arises in the scenario of mixed boundary conditions in which a single absorber featuring a combination of Dirichlet and Neumann has singularities in the normal derivative at their meeting points, even along smooth boundaries. In certain scenarios, the singularity structure is known \cite{BL2018}, hinting that recent methods for rational representation of functions $z^{\alpha}$ may be applicable to the solutions of these cases \cite{Herremans23,Trefethen2021}.

\hl{The approach developed in this paper can also be extended to other rational approximation methods such as AAA-LSQR \cite{cos-tre2023}, the methods of fundamental solutions (see \cite{barnett2008stability} for progress on the Helmholtz equation), or the previously mentioned earlier solution via expansion by a Vandermonde basis orthogonalized by the Arnoldi process. Additionally, it is natural to ask whether a log-lightning method, known to provide exponential convergence in Laplace problems \cite{Trefethen24,Baddoo2021}, can be developed for Helmholtz and ultimately the heat equation.}


Finally, we note that when the solution $u(z,t)$ of \eqref{eq:main} is desired at many $t$ values,  there are certainly economies to be gained by reusing transform evaluations. The current Talbot contour \eqref{eq:Talbot} is robust across a wide range of values $t\in(0,\infty)$, but recomputes the coefficients $\{\hat{u}(z,s_j)\}_{j=1}^M$ at each time point. However, we have observed that other contour choices based on fixed coefficients, can be effective over large time ranges, thus hinting at the potential for an adaptive approach.

\section*{Acknowledgments} AEL acknowledges support from NSF grant DMS-2052636. HL acknowledges support of the Arthur J. Schmitt Presidential Leadership fellowship. The authors are grateful to Prof.~Abinand Gopal for many useful discussions \al{ and two anonymous reviewers for their careful reading and constructive comments on the manuscript}.

\section*{Data availability statement} This work did not involve the use of data.

\section*{Declarations} The authors have no competing interests to declare that are relevant to the content of this article.

\appendix

\section{The logarithmic capacitance}\label{sec:log_cap}

The logarithmic capacitance $d_k$ is determined from the solution of the electrified disk problem
\bsub
\begin{gather}
\Delta v_k = 0; \quad \mbox{in} \quad \mathbb{C}\setminus \mathcal{A}_k; \qquad v_k = 0, \quad \mbox{on} \quad \partial \mathcal{A}_k;\\[4pt]
v_k = \log |z| -\log d_k + \mathcal{O}(|z|^{-1}), \qquad |z|\to\infty.
\end{gather}
\esub
For a variety of regular shapes, the value of $d_k$ is known explicitly. Some examples include (see also \cite{Venu2015})
\bsub
\begin{align}
 d_k &= \frac{\Gamma(\frac12) }{4\pi^{\frac32}}h\approx 0.590h& \mbox{Square of side-length } h.\\
 d_k &= \frac{\sqrt{3}\Gamma(\frac13)^3}{8\pi^2}h\approx 0.422h& \mbox{Equilateral triangle of side-length } h.\\
 d_k &= \frac{3^{\frac34}\Gamma(\frac14)^2}{2^{\frac72}\pi^{\frac32}}h\approx 0.476h& \mbox{Isosceles right triangle with small side-length } h.
\end{align}
\esub
When the shape of the body is completely general, there are several numerical algorithms that can approximate the capacitance \cite{Liesen2017,Baddoo2021}. 

\bibliographystyle{spmpsci}
\bibliography{refs.bib}{}

\end{document}